\documentclass[3p]{elsarticle1}

\newproof{pf}{\bf Proof}

\usepackage{amssymb}
\usepackage{amsmath}
\usepackage{amsthm}

\theoremstyle{plain}
\newtheorem{prop}{Proposition}[section]
\newtheorem{lem}[prop]{Lemma}
\newtheorem{thm}[prop]{Theorem}
\newtheorem{cor}[prop]{Corollary}
\newtheorem{rem}[prop]{Remark}

\theoremstyle{definition}

\newtheorem{exm}[prop]{Example}

\numberwithin{equation}{section}

\newcommand{\exb}{\mbox{\sf Exc}({\mathcal U}_{\beta})}

\newcommand{\potb}{\mbox{\sf Pot}({\mathcal U}_\beta)}

\newcommand{\R}{\mathbb{R}}
\newcommand{\Q}{\mathbb{Q}}
\newcommand{\N}{\mathbb{N}}

\newcommand{\cb}{{\mathcal B}}

\newcommand{\ce}{{\mathcal E}}
\newcommand{\cf}{{\mathcal F}}
\newcommand{\cl}{{\mathcal L}}
\newcommand{\cu}{{\mathcal U}}
\newcommand{\cw}{{\mathcal W}}

\newcommand{\ds}{\displaystyle}

\newcommand{\mcB}{\mathcal B}
\newcommand*{\pairing}[3]{\sideset{_{#1}^{}}{_{#3}^{}}{\mathop{\bracket{{#2}}}}}
\newcommand{\bracket}[1]{\langle #1\rangle}
\newcommand{\norm}[1]{\left\| #1 \right\|}
\renewcommand{\l}{\lambda}
\newcommand{\C}{\mathbb{C}}
\renewcommand{\b}{\beta}

\newcommand{\ve}{\varepsilon}

\begin{document}

\begin{frontmatter}



\title{Potential theory of infinite dimensional L\'evy \mbox{processes}}


\author{Lucian Beznea}
\address{"Simion Stoilow" Institute of Mathematics of
the Romanian Academy, P.O. Box \mbox{1-764,}  RO-014700 Bucharest,
Romania}
\ead{lucian.beznea@imar.ro}

\author{\framebox{Aurel Cornea}\footnote{Aurel Cornea passed away tragically
on September  3$^{\mbox{rd}}$, 2005.}}
\address{Katholische Universit\"at Eichst\"att-Ingolstadt, D-85071 Eichst\"att, Germany}

\author{Michael R\"ockner}
\address{Fakult\"at f\"ur Mathematik, Universit\"at Bielefeld,
Postfach 100 131, D-33501 Bielefeld, Germany, and
Departments of Mathematics and Statistics,
Purdue University,
150 N. University St. West Lafayette, IN 47907-2067, USA.}
\ead{roeckner@mathematik.uni-bielefeld.de}

\begin{abstract}
We study the potential theory of a large class of
infinite dimensional L\'evy processes, including Brownian motion
on abstract Wiener spaces. The key result is the construction of
compact Lyapunov functions, i.e. excessive functions with compact
level sets. Then many techniques from classical potential theory
carry over to this infinite dimensional setting. Thus a number of
potential theoretic properties and principles can be proved,
answering long standing open problems even for the Brownian motion
on abstract Wiener space, as e.g. formulated by R. Carmona in
1980. In particular, we prove the analog of the known result, that
the Cameron-Martin space is polar, in the L\'evy case and apply
the technique of controlled convergence to solve the Dirichlet problem
with general (not necessarily continuous) boundary data.

\end{abstract}


\end{frontmatter}


\section{Introduction}

The purpose of this paper is to study the potential theory of
infinite dimensional L\'evy processes. Such processes, in
particular, the special case of infinite dimensional Brownian
motion, are of fundamental importance as driving (i.e. noise)
processes for stochastic partial differential equations. In
addition, there had been interest in solving Dirichlet problems
for infinite dimensional Ornstein-Uhlenbeck processes (see
\cite{DaGo97}). Nevertheless, there are very few papers in the
last 30 years analyzing these fundamental processes in infinte
dimensions from a potential theoretic point of view, as e.g. in
the nice papers \cite{PrZa04} and \cite{PrZa06} 
on Liouville properties for the
Ornstein-Uhlenbeck process with L\'evy noise. Therefore, many
questions about the validity of fundamental potential theoretic
properties and principles even in the case of Brownian motion on
abstract Wiener space remained open problems, since they were
posed e.g. in \cite{Ca80}, and the more so for infinite
dimensional L\'evy processes.

In this paper we shall establish a number of such properties and
principles answering positively a substantial number of R.
Carmona's questions in \cite{Ca80}. Naturally, in the meantime the
"technology" and methodology in potential theory, in particular,
in its analytic component, has been developed much further (see
e.g. \cite{BeBo04}). The main tool, however, to make this modern
analytic potential theory work in our situation, is the
construction of explicit compact Lyapunov functions, i.e.
($\beta$-) excessive functions with compact level sets, which is
done in a very explicit way for the first time in this paper.
Through such functions the usual local compactness assumption on
the topology can be avoided.

The structure and main results of this paper are the following:

In Section 2 we start with the case of Brownian motion on abstract Wiener
space. The compact Lyapunov functions are constructed in Proposition \ref{prop:p6} and
Theorem \ref{thm:t8}. First consequences are presented in Theorem \ref{thm:2.7} and
Remark \ref{rem.2.9}.
The crucial integrability of the norm $q_x$ (cf. $(2.6)$) with respect to the
Gaussian measure follows from an application of Fernique's Theorem (see
Proposition \ref{prop:p6} ($iv$)).

Section 3 is devoted to infinite dimensional L\'evy processes. The
explicit compact Lyapunov functions are constructed in Proposition
\ref{prop.3.3} and Theorem \ref{th3.5}. Because of lack of an
analog of Fernique's Theorem in this case, we can only consider
Hilbert state spaces and require the existence of weak second
moments (see assumption (H)(i) in Section 3 below). Examples include
perturbations of nondegenerate Gaussian cases and the Poisson case
(see Examples \ref{3.1} and \ref{exm3.6}).

In Section 4 we present the potential theoretic consequences. We
here mention the most important ones only: ($a$) we prove that
Meyer's Hypothesis ($L$) (i.e. existence of a reference measure
for the resolvent) does not hold; ($b$) we derive a natural
condition ensuring that points are polar; ($c$) we prove that the
"Cameron-Martin space" $H$ is polar (including the L\'evy case);
($d$) we introduce natural Choquet capacities (replacing the
Newton capacity in finite dimensions) and show their tightness;
($e$) we prove quasi continuity properties for the excessive
functions; ($f$) we prove the existence of bounded functions
invariant under the semigroup; ($g$) we prove that the state space
$E$ can be decomposed into an uncountable union of disjoint affine
spaces each being invariant under the L\'evy process (Brownian
motion respectively) and that the restriction of the process to
any of such affine subspace is c{\`a}dl{\`a}g; ($h$) we prove that
the so-called "balayage principle" holds.

Results ($d$) and ($h$) above are even new in the infinite
dimensional Brownian motion case.

Section 5 is devoted to the so-called "controlled convergence" for the solution to the
Dirichlet problem for strongly regular open subsets  of $E$.
This type of convergence provides a way to describe the boundary behavior of the solution
to the Dirichlet problem for general (not necessarily continuous) boundary data.
Our main result here is Theorem \ref{thm:t9}.

Finally, we would like to point out that many of the above potential theoretic results
extend to infinite dimensional $\alpha$-stable or more general processes obtained by the above
ones by standard subordination.
In particular, if one considers processes subordinate to infinite dimensional Brownian
motion, such as $\alpha$-stable processes, one can cover jump processes without any conditions on
their weak moments.
We thank Masha Gordina and Sergio Albeverio for pointing this out to us.
More details on this will be the subject of forthcoming work.

In the Appendix we prove a type of analogue to the necessity-part of L. Gross
famous result on measurable norms (see \cite{Gr65}) in the non-Gaussian case.

\section{Brownian motion on abstract Wiener space}

Let $(E,H,\mu)$ be an \emph{abstract Wiener space}, i.e.\
$\bigl(H, \langle\;,\;\rangle\bigr)$ is a separable real Hilbert
space with corresponding norm $\lvert\,\cdot\,\rvert$, which is
continuously and densely embedded into a Banach space
$\bigl(E,\lVert\,\cdot\,\rVert\bigr)$, which is hence also
separable; $\mu$ is a Gaussian measure on $\mathcal B$ ($=$ the
Borel $\sigma$-algebra of $E$), that is, each $l\in E'$, the dual
space of $E$, is normally distributed with mean zero and variance
$\lvert l \rvert^2$. Here we use the standard continuous and dense
embeddings
$$
  E' \subset (H' \equiv) H \subset E\;.
$$
Clearly, we then have that
$$
\leqno{(2.1)}\quad\quad\quad\quad\qquad  {} _{E'}\langle l, h \rangle_E=
\langle l, h \rangle \mbox{ for all }
l\in E' \mbox{ and } h\in H.
$$

We recall that the embedding $H\subset E$ is automatically compact
(see Ch.III, Section~2 in \cite{Bo98}) and that $\mu$ is
$H$-quasi-invariant, that is for $T_h(z):= z+h$, $z,h\in E$, we
have
$$
\mu\circ T_h^{-1} \ll \mu \quad \mbox{ for all }  h\in H.
$$
By the famous Dudley-Feldman-Le Cam Theorem (see \cite{DFLC71} and also Theorem 4.1
in \cite{Ya89} for a concise presentation) we know that the norm
$\lVert\,\cdot\,\rVert$ is  $\mu$-\emph{measurable} in the sense of
L.~Gross (cf.\ \cite{Gr67},  see also \cite{Kuo75}).
Hence also the centered Gaussian
measures $\mu_t$, $t>0$, exist on $\mathcal B$, whose variance are
given by $t\lvert l\rvert^2$, $l\in E'$, $t>0$. So,
$$
\mu_1=\mu\;.
$$
Clearly, $\mu_t$ is the image measure of $\mu$ under
the map $z\longmapsto \sqrt t z$, $z\in E$.

For $x\in E$, the probability measure $p_t(x,\,\cdot\,)$ is defined
by
$$
  p_t(x,A) := \mu_t(A-x)\;
  \quad \text{for all $A\in \mathcal B$.}
$$
Let $(P_t)_{t>0}$ be the associated family of Markovian kernels:
$$
  P_t f(x)
  := \int_E f(y) \; p_t(x,\mathrm dy)
  = \int_E f(x+y) \; \mu_t(\mathrm dy)\;,
  \quad f\in p{\mathcal B}, \; x\in E;
$$
we have denoted by $p\cb$
the set of all positive, numerical, $\cb$-measurable functions on $E$. By Proposition~6 in \cite{Gr67} it follows that $(P_t)_{t\geq 0}$
(where $P_0:= Id_E$) induces a strongly continuous semigroup of
contractions on the space $\mathcal C_u(E)$ of all bounded
uniformly continuous real-valued functions on $E$.

Let ${\mathcal U}=(U_\alpha)_{\alpha>0}$ be the associated Markovian resolvent
of kernels on $(E, \cb)$ given by
$U_{\alpha}:=\int^{\infty}_{0}e^{-\alpha t} P_t\, {\mathrm d}t$, $\alpha>0$.
Recall that $\cu=(U_{\alpha})_{\alpha>0}$ induces a
strongly continuous resolvent of contractions on $\mathcal C_u(E)$.
 By $\mathcal E(\mathcal U)$ we denote the set of all
$\mathcal B$-measurable \emph{$\mathcal U$-excessive functions:}
$u\in\mathcal E({\mathcal U})$ if and only if $u$ is a positive
numerical $\mathcal B$-measurable function, $\alpha U_\alpha u\leq
u$ for all $\alpha>0$ and $\lim_{\alpha\to\infty} \alpha U_\alpha
u(x)=u(x)$ for all $x\in E$.
By Remark~3.5 in \cite{Gr67} it follows that the potential kernel $U$ defined by
$$
 Uf = \int_0^\infty P_tf\;\mathrm dt
$$
is proper, that is, there exists a bounded strictly positive
$\mathcal B$-measurable function $f$ such that $Uf$ is finite.

If $\beta>0$ we denote by $\mathcal
U_\beta$ the sub-Markovian resolvent of kernels
$(U_{\beta+\alpha})_{\alpha>0}$.
Our first aim is to construct a $\cu_\beta$-excessive function $v$ such that:
{\it the set $[v \leq \alpha]$ is  relatively compact for all $\alpha >0$}
(and having some further useful properties).
Such a function will be called {\it\bf compact Lyapunov function} further on.\\

Consider an orthonormal basis
$\{e_n : n\in \N\}$  of $H$ in $E{'}$
which separates the points of $E$.
For each $n\in {\N}$ define
$\widetilde{P}_n :E\longrightarrow
H_n:=span\{e_1,\ldots,e_n\}\subset E{'}$ by

$$
\leqno{(2.2)}\quad\quad\quad\quad\qquad\ds \widetilde{P}_n z=\sum^{n}_{k=1}\ \
{}_{E{'}}\langle e_k,z \rangle_{E}\, e_k,\; z \in E,
$$
and $P_n:=\widetilde{P}_n \upharpoonright_{H}$, so
$$
{P_n}h=\sum^{n}_{k=1}\ \langle e_k,h \rangle \, e_k, \ \
\ h \in H,
$$
and $P_n \longrightarrow Id_H$ strongly as $n\to \infty$.

\begin{lem} \label{lem.2.1}  
$(i)$  Let $y,z\in E.$  Then
$$
{}_{E{'}}\langle \widetilde{P}_n z,y \rangle_{E} =
\sum_{k=1}^n
{}_{E{'}}\langle e_k , z \rangle_{E }\, {}_{E{'}}\langle e_k ,y \rangle_{E}
= {}_{E{'}}\langle \widetilde{P}_n y,z \rangle_{E}.
$$
$(ii)$ Let $y\in E'$, $z\in E$. Then
$$
{}_{E{'}}\langle \widetilde{P}_n y,z \rangle_{E} =
\langle y, \widetilde{P}_n z\rangle  ={}_{E{'}}\langle y, \widetilde{P}_n z \rangle_{E} .
$$
$(iii)$ For $n\geq m$ we have $\widetilde{P}_n\widetilde{P}_m=\widetilde{P}_m\widetilde{P}_n=\widetilde{P}_m.$
\end{lem}

\begin{pf}
The proof of $(i)$ is elementary and that of $(ii)$ follows from $(i)$  and $(2.1)$.
$(iii)$ in turn is a consequence of $(ii)$.
\hfill $\square$
\end{pf}

\begin{prop}\label{prop:p5} 
We have
$$
\lim_{n\rightarrow\infty}|| \widetilde{P}_{n}z-z ||=0  \mbox { in }
\mu\mbox{-measure.}
$$
\end{prop}

\begin{pf}
Let $\nu$ be the cylinder measure on $H$ corresponding to $\mu$.
Let $i:H\longrightarrow E$ denote the above embedding. Then again
by the Dudley-Feldman-LeCam Theorem (\cite{Ya89}, Theorem $4.1,$ in
particular $(iv)$) for all $\varepsilon>0$
$$
\leqno{(2.3)}\ \ \ \ \ \lim_{m,n\to\infty} \nu (\{{h \in H : \,
||P_nh-P_mh||>\varepsilon} \} )=0.
$$
But $\mu (\{{z \in E :  \, ||\widetilde{P}_nz-\widetilde{P}_m
z||>\varepsilon} \} )= \nu (\{ {h \in H : \,
||P_nh-P_mh||>\varepsilon} \})$.
Hence by $(2.3)$ there exists a
$\cb(E)/\cb(E)$-measurable function $F:E\longrightarrow E$ such
that
\begin{center}
$\ds \lim_{n\rightarrow\infty}||F-\widetilde{P_n}||_{E}=0$ in
$\mu$-measure,
\end{center}
and therefore $\mu$-a.e. for a subsequence $(n_k)_{k \in \N}$.
Thus for all $m \in {\N}$ and $\mu$-a.e. $z\in E$
$$
{}_{E'}\langle e_m,z\rangle_E=\lim_{k\to \infty} {}_{E'}\langle
e_m,\widetilde{P}_{n_k} z \rangle_E= {}_{E'}\langle
e_m,F(z)\rangle_E \, ,
$$
and we conclude that $F(z)=z$ for $\mu$-a.e. $z\in E$. 
\hfill $\square$
\end{pf}

Passing to a subsequence if necessary, which we denote by $Q_n$,
$\widetilde Q_n$,
$n \in {\N}$, respectively, we may assume that
$$
\leqno{(2.4)}\quad\quad \quad ||Id_H-  Q_n||_{\cl(H,E)} \leq \frac{1}{2^n}
$$
and
$$
\leqno{(2.5)}\quad\quad \quad \mu(\{z \in E : \,  ||z-\widetilde{Q}_nz|| >
 \frac{1}{2^n} \}) \leq  \frac{1}{2^n},
$$
where we used the compactness of the embedding $H \subset E$ for
$(2.4)$ and Proposition \ref{prop:p5}  for $(2.5)$.\\

Let $x\in E\setminus H$.
We note that assuming  the existence of such a point implies that $\dim H=\infty$
and a standard argument shows that $\mu(H)=0$ (see \cite{Bo98}).\\

The following lemma is due to R. Carmona.

\begin{lem}\label{lem:3.2} 
Let $x \in E \setminus H$.
There exists an orthonormal basis
$\{e_n^{x}\, :\, n\in \N\}$ of $H$ such that $e_n^{x}\in E^{'}$ for
all $n \in \N$,  $\{e_n^{x}\, : \, n\in \N\}$ separates the points of
$E$ and
$$
{_{E{'}}} \langle e^x_n,x \rangle_{E}\geq 2^{\frac{n}{2}} \mbox{ for all } n.
$$
\end{lem}

\begin{pf}
This follows from the proof of Lemma 1 in {\rm \cite{Ca80}}.
Concerning the claim
that $\{e_n^{x}\, :\, n\in \N\}$ separates the points of $E$, one just realizes
that $\{x_n^{*} \, : \, n\in \N\}$ in the proof of
 \cite{Ca80}, Proposition 1,  separates the points of $E$ and it
follows by the construction there, that so does $\{e_n^x\, :\,
n\in\N\}$.
\hfill $\square$
\end{pf}

Define the function $q_x :E\longrightarrow \overline{{\R}}_{+}$ by
$$
\leqno{(2.6)}\quad\  \ds q_x (z):=  \left[ \sum_{n\geq 0} 2^n ||\widetilde{Q}_{n+1} z - \widetilde {Q}_{n}z ||^2 +
\left( \sum_{n\geq 1}  2^{-\frac{n}{2}}   | {}_{E'}\langle e^x_n,  z \rangle_{E} |  \right)^2 \right]^{\frac{1}{2}},\quad
z\in E,
$$
where $\widetilde{Q}_0:=0$ and $e_n^{x},$ $n\in \N,$ is as defined in Lemma \ref{lem:3.2}.
Also $\widetilde{Q_n}$, $n\in\N$, is defined as above with this particular ONB.
Let
$$
\ds E_x:=\{z\in E :    q_x (z)<\infty\}.
$$
Note that by Lemma \ref{lem:3.2} we have
$$
x\in E\setminus E_x.
$$

Recall that if $l \in E{'}$ then for all $z \in E$ we have
$$
\leqno{(2.7)}\quad\quad\quad  \displaystyle\int_{E}{l^2(y)\, p_t(z,dy)}=t|l|^2+l^2(z),
$$
where $|l|$ denotes the $H$-norm of $l$ ($\in E'\subset H' \equiv H$).\\

Modifying the arguments in \cite{Ku82} we can now prove:

\begin{prop}\label{prop:p6}  
Let $x\in E \setminus H$. The following assertions  hold.\\
$(i)$  $\mu(E_x)=1$.\\
$(ii)$  For all $h \in H$ we have $q_x (h)\leq \sqrt{3} |h|$. In
particular, $H \subset E_x$
continuously.\\
$(iii)$  For all $z \in E$ we have
$$
||z|| \leq\sqrt{2} \, q_x (z).
$$
In particular, $(E_x,  q_x )$ is complete. Furthermore, $(E_x, q_x
)$
is compactly embedded into $(E, ||\cdot||)$.\\
$(iv)$ $(E_{x},H,\mu)$ is an abstract Wiener space. In particular,
$q_x\in L^2(E,\mu)$.
\end{prop}

\begin{pf}
$(i)$ Let us set
$$
g(z) :=  \sum_{n\geq1} 2^{-\frac{n}{2}}   | {}_{E'}\langle e^x_n,  z \rangle_{E} | , \quad z\in E.
$$
We show that
$$
\leqno{(2.8)}\quad\quad g\in L^2(E,\mu).
$$
Indeed, by $(2.7)$ and Minkowski's inequality we have

$$
\int_E g^2(z)\mu(dz)\leq
\left( \sum_{n\geq1} 2^{-\frac{n}{2}} \sqrt{\int_E {}_{E'}\langle e^x_n,  z \rangle_{E}^2}\,  \mu(dz) \right)^2 =
\left(\sum_{n\geq1} 2^{-\frac{n}{2}} |e^x_n|\right)^2 < \infty.
$$
Consequently $g$ is finite $\mu$-a.s. and  assertion $(i)$ is now
a direct consequence of $(2.5)$ and the Borel-Cantelli Lemma.

$(ii)$ For all $h \in H$, by $(2.4)$, we have
$$
\ds ||\widetilde{Q}_{n+1}h-\widetilde{Q}_{n}h|| \leq
2^{-n}|Q_{n+1}h| \leq 2^{-n} |h|
$$
and therefore
$$
q_x (h)^2 \leq  \sum_{n\geq 0} 2^{-n} |h|^2 +(\sum_{n=1}^\infty 2^{-n})
\sum_{n=1}^\infty  \langle e^x_n , h\rangle^2,
$$
which implies the  assertions of $(ii)$.

$(iii)$  We have for all $n \in {\N}$ and $z \in E$
$$
\sup_{m\geq n}||\widetilde{Q}_mz-\widetilde{Q}_nz|| \leq
\sup_{m\geq n}\sum^{m-1}_{k=n}||\widetilde{Q}_{k+1}z-\widetilde{Q}_{k}z|| 2^{\frac{k}{2}} 2^{-\frac{k}{2}}
\leq
$$
$$
(\sum^{\infty}_{k=n} 2^{k}||\widetilde{Q}_{k+1}z-\widetilde{Q}_{k}z||^2)^{\frac{1}{2}}
(\sum^{\infty}_{k=n} 2^{-k})^{\frac{1}{2}} \leq
q_x (z)(\sum^{\infty}_{k=n} 2^{-k})^{\frac{1}{2}}.
$$
In particular (restricting the above to $z \in E_x$),
$(\widetilde{Q}_n)_{n \in {\N}}$ is a Cauchy sequence in
$\mathcal{L}(E_x, E)$ with respect to the operator norm. Hence by
completeness there exists $T \in \mathcal{L}(E_x , E)$ such that
$\widetilde{Q}_n\rightarrow T$ as $n\to \infty$ in operator norm
and $T$ is compact since each $\widetilde{Q}_n$ is of finite rank.
By Lemma \ref{lem.2.1} $(ii)$ it follows that for each $e^x_n$
$$
{}_{E{'}}\langle e^x_n, Tz \rangle_{E}=
\lim_{m\to \infty}{}_{E{'}}\langle e^x_n, \widetilde Q_m z \rangle_{E}=
\lim_{m\to \infty}{}_{E{'}}\langle Q_m e^x_n,  z \rangle_{E}=
{}_{E{'}}\langle e^x_n, z \rangle_{E}.
$$
Therefore, for all $z\in E_x$, $Tz=z$ and thus $E_x\subset E$
compactly and furthermore
$$
||z||=||Tz||=\lim_m||\widetilde{Q}_mz||=\lim_m ||\widetilde{Q}_m z-\widetilde{Q}_0 z||\leq
\sup_{m\geq 0}||\widetilde{Q}_m
z-\widetilde{Q}_0 z|| \leq
q_x (z)(\sum^{\infty}_{k=0} 2^{-k})^{\frac{1}{2}}.
$$
The completeness of $(E_x, q_x)$ then easily  follows by Fatou's
lemma.

$(iv)$ {\it Claim 1.} Let $z\in E_x$. Then $\ds
\lim_{n\rightarrow\infty} q_x(z-\widetilde{Q_n}z)=0$. In
particular, $H\subset E_{x}$ densely.

{\it Proof of Claim 1.} For all $n\in\N$ by Lemma \ref{lem.2.1} $(ii)$ and $(iii)$
$$
q_x^2(z-\widetilde{Q}_n z)=
\sum_{k=0}^{\infty}2^k ||\widetilde{Q}_{k+1}z-\widetilde{Q}_{(k+1)\land n}z-
\widetilde{Q}_k z + \widetilde{Q}_{k \land n} z||^2+
\sum_{k=1}^{\infty} 2^{-\frac{k}{2}}
|{}_{E'}\langle (Id_H - Q_n)e_k^x,z  \rangle_{E}|
$$
$$
=\sum_{k\geq n} 2^k ||\widetilde{Q}_{k+1}z-\widetilde{Q}_k z||^2+
\sum_{k \geq N_n} 2^{-\frac{k}{2}}|{}_{E'}\langle e_k^x,z  \rangle_{E}|
$$
for some $N_n\nearrow\infty$ when $n\rightarrow\infty$. Now the
first part of the assertion follows, since $z\in E_{x}$. The
second part is then a consequence thereof, since $\widetilde{Q}_n
z\in H$ for all $n\in\N$.

{\it Claim 2.} Let $l\in E'_x$ and $l_n:=l \circ \widetilde{Q}_n$,
$n\in\N$. Then $l_n\in E'$ and $\ds \lim_{n\rightarrow\infty}
l_n(z)=l(z)$ for all $z\in E_x$.

{\it Proof of Claim 2.} Since each $\widetilde{Q}_n:E\rightarrow
H$ is continuous and $H \subset E_{x}$ continuously, we have that
$l_n\in E'$ for all $n\in\N$. The last part of the assertion
follows from Claim 1.\\

We shall now see that Claim 1 and Claim 2 imply assertion $(iv)$.
Indeed, since $H\subset E_x$ continuously by $(ii)$ and densely by
Claim 1, it remains to show that $\mu$ is centered Gaussian as a
measure on the Banach space $(E_x,q_x)$, with Cameron-Martin space
$H$, i.e. every $l\in E'_x$ has a mean zero normal distribution with
variance $|l|^2$. (Recall that $E'_x\subset (H{'}\equiv)H\subset
E_x$ continuously and densely.)  So, let $l\in  E'_x$ and let
$l_n$, $n\in\N$, be as in Claim 2. Then $l_n$, $n\in\N$, are jointly
Gaussian with mean zero and $l_n\longrightarrow l$ $\mu$-a.e. as
$n\longrightarrow\infty$ by $(i)$, hence $l_n\longrightarrow l$ in
$L^2(E,\mu)$ as $n\longrightarrow\infty$. Since then $l_n
\longrightarrow h$ in $H$ as $n\longrightarrow\infty$ for some
$h\in H$, considering the Fourier transforms we see that $l$ under
$\mu$ has a mean zero normal distribution with variance $|h|^2$.
But obviously $l_n \longrightarrow l$ weakly in $H$, hence $l=h$.
The last part of assertion $(iii)$ then follows by Fernique's
Theorem (see e.g. [10, Theorem 2.8.5]).
\hfill $\square$

\end{pf}

\begin{cor} \label{cor.2.4}   
{\rm (cf. \cite{Ca80}, Proposition 1)}  Let  $x\in E\setminus H$.
Then  there exists a Borel linear subspace  $E_x $ of $E$  such
that $H\subset E_x$, $\mu(E_x)=1$, and $x\not\in E_x$. In
particular, $\mu(H+x)=0$
\end{cor}

\begin{pf}
The first part is just Proposition \ref{prop:p6} $(i)$. Since
$(x+H)\cap E_x=\emptyset$, also the second part of the assertion
follows.
\hfill $\square$
\end{pf}

\begin{lem}\label{lem:l7}  
Let $L \in \cb$ be a linear subspace of $E$ such that $\mu(L)=1$.
Then for all $z \in E$ the set $L+z$ is invariant with respect to
$(P_t)_{t\geq{0}}$, i.e. $P_t(1_{L+z})=1_{L+z}$ for all $t>0$.
In particular, the measure $p_t(x,\cdot)$
is carried by $L+z$ for every $x \in{L+z}$.
\end{lem}

\begin{pf}
We have $\mu_{t}(L)=\mu_1(t^{-\frac{1}{2}}L)=\mu(L)=1$. Let $z \in
E$. If $x \in{L+z}$ then $p_t(x,L+z)=\mu_t{(L+z-x)}=\mu_{t}(L)=1$.
If $x \notin{L+z}$ then $(L+z-x)\cap{L}=\emptyset$ and thus
$p_t(x,L+z)=\mu_t{(L+z-x)}\leq \mu_t{(E \setminus L)}=0$.

\hfill $\square$
\end{pf}

\begin{thm}\label{thm:t8}  
Let $x\in E \setminus H$. Define $v_0^x:=U_1 q^2_x $ and for every
$z\in E$, $v_z^x:=v_0^x\circ T_z^{-1}$. Then $v_z^x$ is a compact
Lyapunov function such that  $E_x+z=[v_z^x< \infty]$ and each
$E_x+z$ is invariant with respect to $(P_t)_{t\geq 0}$.
\end{thm}

\begin{pf}
By Proposition \ref{prop:p6} and Lemma \ref{lem:l7} it follows
that  $E_x+{z}$ is absorbing and invariant with respect to
$(P_t)_{t\geq 0}$.

We show that $v_0^x$  is a compact Lyapunov function on $E$ such
that $E_x =[v_0^x<\infty]$. By Proposition \ref{prop:p6} $(iv)$
and by $(2.8)$ we have $q_x \in L^2(E, \mu)$. Let $\ds M:=\int_{E}
q^2_x (y)\mu(dy)$. Then for all $t>0$, $z \in E$, $\ds
\int_{E}q^2_x (y)\mu_t(dy)= {M}{t}$, and by the sublinearity of
$q_x$
$$
P_t(q^2_x )(z)=\int_{E}q^2_x (z+y)\mu_{t}(dy)\leq
2\int_{E} (q^2_x (z)+q_x^2(y))\mu_{t}(dy)\leq
 2(q^2_x (z) +Mt).
$$
We conclude that
$$
v_0^x(z)=U_1 (q^2_x )(z)= \int^{\infty}_{0}e^{-t}P_t(q^2_x )(z)dt\leq
2q^2_x (z)+2M\int_0^\infty e^{-t}t dt.
$$
Hence $E_x \subset [v_0^x< \infty]$.

We claim that $v_0$ has compact level sets in $E$.
Obviously, $q_x$ is lower semicontinuous on $E$.
Therefore, because $U_1$ maps bounded continuous functions to bounded continuous functions,
$v_0^x$ is also lower semicontinuous on $E$.
Then by Proposition \ref{prop:p6} the sets $[q_x \leq \beta]$ are compact in $E$,
hence it will be sufficient to prove that
$$
v_0^x \geq q^2_x .
$$
Let
$f_n(z):=||\widetilde{Q}_{n+1}z-\widetilde{Q}_{n}z||^2$ and
$(l_k)_k \subset E{'}$, $||l_k||=1$, be such that for all $z \in
E$
\begin{center}
$\ds
||z||=\sup_{k}l_k(z)$.
\end{center}
The functionals $l_{k,n}:=l_k \circ
(\widetilde{Q}_{n+1}-\widetilde{Q}_{n})$ belong to $E'$ and using
$(2.7)$ we get for all $z\in E$, $t>0$ and natural number $n$:
$$
P_tf_n(z)=\int_{E}f_n(y)p_t(z,dy)=
\int_{E}\sup_{k}l^{2}_{k,n}(y) p_t(z,dy)\geq
$$
$$
\sup_{k}\int_{E}l^{2}_{k,n}(y) p_t(z,dy)\geq
\sup_{k}l^2_{k,n}(z)=f_n(z).
$$
Hence $P_tf_n\geq f_n$. Recall that $g$ denotes the second sum
occurring in the definition of $q_x$. We have
$$
P_t(g^2)(z)\geq (P_t g(z))^2=
\left( \sum_{n\geq 1}  \frac{1}{\,\,\,  2^\frac{n}{2}}
\int_E |{}_{E'}\langle e^x_n,  y \rangle_{E}| \, p_t(z, dy)  \right)^2 \geq
$$
$$
\left( \sum_{n\geq 1}  \frac{1}{\,\,\,  2^\frac{n}{2}}   |
\int_E {}_{E'}\langle e^x_n,  z+y \rangle_{E} \mu_t(dy)| \right)^2=
\left( \sum_{n\geq 1}  \frac{1}{\,\,\,  2^\frac{n}{2}}   |
{}_{E'}\langle e^x_n,  z\rangle_{E}| \right)^2= g^2(z).
$$
Hence  we also have $P_t(g^2)\geq g^2$.
Since $q^2_x = \sum_{n\geq 0} 2^n f_n + g^2$ we obtain
$$
P_t (q^2_x ) \geq q^2_x  \mbox{ for all } t>0
$$
and thus
$$
v_0^x=\int^{\infty}_{0}e^{- t}P_t (q^2_x )\, dt
\geq q^2_x \int^{\infty}_{0} e^{-t}dt=q^2_x  .
$$
Since $P_t(f\circ T_z)=P_tf \circ T_z$ for all $f \in p\cb$ and
$z\in E$, we deduce that if $u \in \ce(\cu_\beta)$ then $u \circ
T_z \in \ce(\cu_{\beta})$.
Consequently, by the first part of the
proof, the function $v_z^x=v_0^x \circ T_{-z}$ is a compact
Lyapunov function for every $z \in E$ and $E_x+{z}=[v^x_z<\infty]$.
\hfill $\square$
\end{pf}

\begin{rem}\label{rem2.7'} 
Fix $x\in E$ and for $y,z\in E$ define the equivalence relation $y
\sim z$ if and only if  $y-z\in E_x$, and let $\tau$ be defined as
a set in $E$ containing exactly one representative of each
equivalence class. Note that since $\alpha x+E_x$, $\alpha\in\R$,
are pairwise disjoint, $\tau$ is uncountable, and
$$
E=\bigcup_{z\in\tau}^\cdot (E_z+x).
$$
Hence $E$ is an uncountable union of disjoint  Borel sets which
are invariant for the Brownian motion.
\end{rem}

As one consequence of Theorem \ref{thm:t8}, we can reprove Gross's famous result
on the existence of the infinite dimensional Brownian motion 
(cf. \cite{Gr67}; see also \cite{Pi71} and \cite{Pi72} for constructions of diffusion processes on abstract Wiener spaces)
and give some additional information,
based on a general technique we developed in \cite{BeRo09}; the proof will be sketched.

Recall that a \emph{Ray cone} associated with $\mathcal U_\beta$, $\beta>0$,
is a cone $\mathcal R$ of bounded $\mathcal U_\beta$-excessive
functions such that:
$U_\alpha(\mathcal R)\subset \mathcal R$ for all $\alpha>0$,
$U_\beta\bigl((\mathcal R-\mathcal R)_+\bigr)\subset \mathcal R$,
$\sigma(\mathcal R)=\mathcal B$,
it is $\min$-stable, separable in the supremum norm and $1\in\mathcal R$.
The topology on $E$ generated by a Ray cone is called
\emph{Ray topology}.

\begin{thm}\label{thm:2.7}  
$(i)$ There exists a diffusion process $\mathcal W =
(\Omega,\mathcal F,\mathcal F_t, W_t,\theta_t,P^x)$ with state
space $E$ ({\rm the  Brownian motion on $E$}), having
$(P_t)_{t\geq 0}$ as transition function.

$(ii)$ The topology of $E$ is a Ray one. For every finite measure  $\lambda$ on $(E, \cb)$ there
exists a natural capacity associated with the Brownian motion on an abstract Wiener
space, which in particular is tight. More precisely,
 the functional $M\longmapsto c_\lambda(M)$, $M\subset E$, defined by
$$
c_\lambda(M) := \inf \bigl\{ \lambda(P_{T_G} p) :  M\subset G\;\text{open} \bigr\}
$$
is a Choquet capacity on $E$, where
$P_{T_G}$ denotes the hitting kernel of the set $G$
(see, e.g.,  Section 5 below for further details)
and $p$ is a bounded $\cu$-excessive function of the form
$p=Uf_0$ with $f_0\in bp\cb$  strictly positive;
$bp\cb$ denotes the bounded elements of $p\cb$.

$(iii)$  Every $\cu$-excessive function $u$ of the form $u=Uf$, $f\in p\cb$,  is
$c_\lambda$-quasi continuous, provided it is  finite $\lambda$-a.e.
More generally, every potential of a continuous additive functional
(cf. \cite{Sh88} or \cite{BeBo04})
is
$c_\lambda$-quasi continuous if it is finite $\lambda$-a.e.
In particular, every $\cu$-excessive function is $c_\lambda$-quasi lower semicontinuous.
\end{thm}

\noindent
{\bf Sketch of the proof.}
$(i)$ We show first that $\cu$ satisfies condition $(*)$ from \cite{BeRo09}, Corollary 5.4, namely
for some $\beta >0$ and every $z\in E$ we have:

$(*)\quad$ {\rm if $\xi \in \exb$ and
$\xi \leq U_{\beta}(z, \cdot)$ then $\xi \in \potb$};\\
we have denoted by $\exb$ (resp. $\potb$) the set of all
$\cu_\beta$-excessive measures (resp. of all potential
$\cu_\beta$-excessive measures). Let $x, z\in E$. Theorem
\ref{thm:t8} and  assertion $(ii)$ of Corollary 5.4  from \cite
{BeRo09} imply that the restriction of $\cu$ to $E_x+{z}$ is the
resolvent of a right process with state space $E_x+{z}$. Therefore
it verifies in particular $(*)$ for $z\in E_x+z$; cf. assertion
$(ii.1)$ of Corollary 5.4   from \cite{BeRo09}. Hence $(*)$ holds
for all $z\in E$ and so, by assertion $(i)$ of Corollary 5.4  in
\cite{BeRo09}, we conclude  now that $(P_t)_{t\geq{0}}$ is the
transition function of a Borel right process with state space $E$.

The argument in \cite{Gr67}, page 134, ensures (using a criterion
of E. Nelson, \cite{Ne58}) that the process has continuous paths.

$(ii)$ Since the semigroup $(P_t)_{t\geq 0}$ is strongly continuous on ${\cal C}_u(E)$,
we deduce from Proposition 2.2 in \cite{BeRo09} that the topology of $E$ is a Ray one.
By the above considerations and Proposition 4.1 in \cite{BeRo09}
we get the desired capacity and its tightness property.

Assertion $(iii)$ is a consequence of Proposition 3.2.6 from \cite{BeBo04},
using essentially the property of the topology to be a Ray one, proved above. \\

\noindent
\begin{rem} \label{rem.2.9}    
$(i)$ The existence of the compact Lyapunov function $v_z^x$ was
crucial in our approach. To underline this, we present here the
main arguments from the proof of Theorem 5.2  from \cite{BeRo09},
on which (the above crucially used) Corollary 5.4 is based: The
resolvent $\cu$ is always associated to a Borel right process, but
on a bigger set $E_1$, the so called "entry space". However, if
there exists a nest of Ray compact sets, then the set,
$E_1\setminus E$ is polar and consequently $\cu$  is the resolvent
of the process restricted to $E$  (see,  e.g., Lemma~3.5 in
\cite{BeBoRo06a}). The level sets $[v_z^x\leq n]$, $n\in \N$,
offer precisely the required  nest of Ray compact subsets of
$E_x+z$ and therefore the restriction of $\cu$ to $E_x+z$ is the
resolvent of a Borel right process with state space $E_x+z$, for
all $x,z\in E$.

$(ii)$  In \cite{Ca80}, page 41, R. Carmona asked whether
there is a relevant notion of Newtonian capacity
in  the setting of the infinite dimensional Brownian motion.
The second  assertion of $(ii)$ in  Theorem \ref{thm:2.7} answers this  question;
see also Section 4 below. The quasi continuity properties stated by assertion $(iii)$ of
Theorem \ref{thm:2.7} are exactly analogous to those which hold in the classical case
with respect to the Newtonian capacity.
\end{rem}

\section{L\'evy processes on Hilbert space}  

The purpose of this section is to show that a slight modification
of the construction in the previous section gives rise to explicit
compact Lyapunov functions for L\'evy processes in infinite
dimensions provided they have finite (weak) second moments. For
simplicity we restrict ourselves to the case of Hilbert state
spaces. As in Section 2 we start with a separable real Hilbert
space $(H, \langle , \rangle)$ with corresponding norm $|\cdot|$
and Borel $\sigma$-algebra $\mcB(H)$.

Let $\l : H \longrightarrow \C $ be a continuous
negative definite function such that $\l (0) =0$. Then
by Bochner's Theorem there exists a finitely
additive measure $\nu_t, \, t > 0,$ on $(H , \mcB(H))$
such that for its Fourier transform we have
$$
\widehat  \nu _t ( \xi ) := \int _H e^{i \langle\xi , h \rangle} \nu _t (dh ) =
e ^{- t \l ( \xi )},\quad \xi \in H.
$$

Let $E $ be a Hilbert space such that $H \subset E $ continuously
and densely, with inner product $\langle , \rangle_E$ and norm
$\norm \cdot $. Then, identifying $H$ with its dual $H'$ we have
$$
\leqno{(3.1)}\quad\quad\quad\quad
E'\subset H \subset E
$$
continuously and densely, and $\pairing{E'}{\xi , h}{E} =\langle \xi, h \rangle$,
for all $\xi \in E ', \, h \in H$.

In addition, we assume that
the following assumption holds
$$
\leqno{(HS)}\quad\quad\quad H\subset E\,  \mbox{ is Hilbert-Schmidt.}
$$
(Such a space $E$ always exists.) Then,
since $\widehat\nu_t$ is continuous on $H$,
by the Bochner-Minlos
Theorem (see, e.g.,  [Ya89]) each $\nu_t$ extends to a measure on
$(E, \mcB(E))$, which we denote again by $\nu_t$, such that
$$
\leqno{(3.2)}\quad\quad\quad \widehat  \nu_{t} (\xi )
=\displaystyle\int _E e^{ i {}_{E'}\! \langle \xi, z\rangle_E }
\nu_t (dz) \mbox{ for all } \xi \in E'.
$$
Clearly, $\l$ restricted to $E'$ is Sazonov continuous, i.e.,
continuous with respect to the topology generated by all Hilbert-Schmidt operators on $E'$. Hence by L\'evy's continuity theorem on
Hilbert spaces (see [35, Theorem IV.3.1 and Proposition
VI.1.1]), $\nu _t \to \delta _0$ weakly as $t \to 0$. Here
$\delta_0$ denotes Dirac measure on $(E, \mcB (E))$ concentrated
at $0\in E$. Furthermore, by the L\'evy-Khintchine Theorem on
Hilbert space (see,  e.g., Theorem VI.4.10 in \cite{Pa67})
$$
\leqno{(3.3)}\quad\l (\xi) = -i {\, } _{E'}\langle \xi, b
\rangle_E + \frac 12 {\, } _{E'} \langle \xi , R \xi \rangle_E -
\displaystyle \int _E
 \left( e^{ i {}_{E'}\! \langle \xi, z\rangle_E } -1  -
 \frac{i {\, } _{E'} \langle \xi , z \rangle_E }{1 + \norm z ^2 } \right) M (dz ) ,  \xi \in E',
$$
where $b \in E$, $ R : E '\longrightarrow E$ is linear such that
its composition $R\circ i_R$
with the Riesz isomorphism $i_E : E \to E'$ is a
non-negative symmetric trace class operator on $E$, and $M$ is a L\'evy
measure on $(E, \mcB(E))$, i.e. a positive measure on $(E,
\mcB(E))$ such that
\[M (\{0\})=0,\quad \int _E ( 1 \wedge \norm z ^2 ) M (dz ) < \infty.\]
Defining the probability measures
$$
\leqno{(3.4)}\quad\quad\quad
 p_t (x, A) := \nu_t (A-x), \;\,  t> 0, \; x \in E, A\in \cb(E),
$$
we obtain a semigroup of Markovian kernels $(P_t)_{t\geq 0}$ on
$(E, \mcB(E))$ just like for the Gaussian case in the previous
section. It has been proved in \cite{FuRo00}, that there exists a
conservative Markov process
$X=(\Omega,\cf,\cf_t,X_t,\theta_t,P^{x})$  with transition
function $(P_t)_{t\geq 0}$ which has c\`adl\`ag paths (see Theorem
5.1 in \cite{FuRo00}). $X$ is just an infinite dimensional version
of a classical L\'evy process. Obviously, each $P_t$ maps $C_b(E)$
into $C_b(E)$, hence so does its associated resolvent $U_\b= \int
_0^\infty e ^{- t \beta} P_t \, dt $, $\beta > 0$. In addition,
$P_t f ( z) \to f (z)$ as $t \to 0$, hence $\beta U_\beta f (z)
\to f (z)$ as $\beta\to \infty$ for all $f \in C_b(E)$, $z \in E$.
Hence $X$ is also quasi-left continuous, and thus a standard
process.

By $(HS)$ there exists an orthonormal basis $\{e_n :\ n\in\N \}$ of $H$ contained in $E'$
having the following properties:\\
There exist $\lambda_n\in (0, \infty)$, $n\in \N, $ such that
$$
\sum_{n=1}^\infty \lambda_n < \infty
$$
and $\overline e_n:= \frac{e_n}{\sqrt{\lambda_n}}$, $n\in \N$,
form an orthonormal basis of $H$. Furthermore,
$$
\leqno{(3.5)}\quad\quad\quad \quad
 \lambda_n {} _{E'}\langle e_n, z \rangle_E=\langle e_n, z \rangle_E\quad
\mbox{ for all } n\in \N, z \in E.
$$
In particular, $\{e_n : \, n\in \N\}$ separates the points of $E$.
The construction of  $\{e_n : \, n\in \N\}$ is standard. We refer, e.g.,
to  Proposition 3.5 from \cite{AR88}. For $n\in \N$ define
$\widetilde P_n : E\longrightarrow E'$ by
$$
\widetilde P_n z :=\sum_{k=1}^n  {} _{E'}\langle e_k, z \rangle_E\, e_k, \, z\in E,
$$
and $P_n:= \widetilde P_n\upharpoonright_{H}$.
Since by $(3.5)$ for all $n\in \N$ and $z\in E$
$$
\widetilde P_n z =\sum_{k=1}^n  \langle \overline e_k, z \rangle_E\, \overline e_k,
$$
we have
$$
\leqno{(3.6)}\quad\quad\quad \quad
\lim_{n\to \infty} || \widetilde P_n z - z ||=0 \mbox{ for all } z\in E.
$$

\begin{rem} \label{rem.3.1} 
Let $t>0$ and consider the (non-Gaussian) triple $(E, H, \nu_t)$.
As mentioned at the beginning of Section 2, in the Gaussian case
the Dudley-Feldman-Le Cam Theorem says that $||\cdot ||$ is a
$\mu$-measurable norm in the sense of Gross, which, however, is
not known to be true for our not necessarily Gaussian measure
$\nu_t$. Recall that in \cite{DFLC71} only a weaker notion of
"$\mu$-measurability" was shown and this notion was proved to be
equivalent with Gross's $\mu$-measurability only in the Gaussian
case (see [16, Theorem 3]). $(3.6)$ above, however, provides
a suitable substitute for the special sequence $(P_n)_{n\in \N}$
of projections considered above, whose existence follows from
assumption $(HS)$. It is an interesting question whether this
depends on this special sequence $(P_n)_{n\in \N}$ or, whether
$(3.6)$ is true at least $\nu_t$-a.s. for any sequence of
projections $(P_n)_{n\in \N}$ of the type considered in Section 2,
i.e., whether Proposition \ref{prop:p5} is true for $\nu_t$ or
even more general measures. This question (of independent
interest) is answered in the Appendix. The corresponding
Proposition {A.2} can be considered as a kind of generalization of
the Dudley-Feldman-Le Cam Theorem to non-Gaussian measures under
assumption $(HS)$.
\end{rem}

Now we want to extend the construction of compact Lyapunov
functions from Section 2 to this
case. To this end we have to make the following further assumption $(H)$
below, which as we shall see (cf. Example~\ref{3.1} below),
is always fulfilled if $\l $ is sufficiently regular.

\begin{enumerate}
   \item [(H)(i)] There exists $C> 0$ such that for all $\xi\in E'$
   \[\int_E {}_{E'}\langle \xi , z \rangle_E ^2 \, \nu _t ( d z ) \leq C ( 1 + t ^2 )|\xi|^2 , \quad t > 0.\]
   \item[(ii)] $\nu_t(H)=0$ for all $t>0$.
\end{enumerate}

\begin{exm}\label{3.1} 
$(i)$ If $\l$ is sufficiently regular, by a straightforward computation one deduces from the
 representation in $(3.3)$  that for every $\xi \in E'$
 \begin{align*}
   \int_E {}_{E'}\langle \xi , z \rangle_E ^2 \nu_t ( d z ) &=
- \frac{ d ^2 }{d\ve ^2}e ^{-t \l ( \ve \xi)}\big | _{\ve =0}\\
   & = t ^2 \left( {}_{E'}\langle \xi , b \rangle_E +
\int _E {}_{E'}\langle \xi ,z \rangle_E \frac{\norm z ^2}{1 + \norm z ^2} \, M (dz )\right)^2\\
   &\quad +t \left(  {}_{E'}\langle\xi , R\xi \rangle_E
+ \int_E {}_{E'}\langle \xi, z  \rangle_E ^2 \,  M (dz)\right)
 \end{align*}
where we assume that $\xi $ is such that $\ds \int_E
{}_{E'}\langle \xi, z \rangle_E ^2 \, M  (dz) < \infty$.
Hence assuming that $b\in H,$ $R(E')\subset H$ and
$R: E' \longrightarrow H$ is continuous with respect to the norm $|\cdot |$ on $E'$,
we have that
$(H) (i)$ holds provided $\ds \int_E {}_{E'}\langle \xi
, z\rangle_E ^2 \,  M (d z )<\infty$ for all $\xi$ in $E{'}$,
because then by the uniform boundedness principle $\ds \sup\{\int
_E {}_{E'}\langle \xi , z\rangle_E ^2 \, M (d z ) : \,  | \xi |
\leq 1\}< \infty$.

$(ii)$  Assume that $\lambda$ is such that in $(3.3)$
$R=i_H \circ i_H^* \circ i_E^{-1}$,
where $i_H$ denotes the embedding $H\subset E$
and $i_H^*: E\longrightarrow H$ its adjoint.
Fix $t>0$.
Then there exist  probability measures
$\mu_t$,  $\nu_t^{0}$ on $(E,\cb(E))$ and $b\in E$ such that
$$
\nu_t=\delta_{tb}* \mu_t * \nu_t^{0},
$$
where $\mu_t$ is Gaussian such that  $(E,H,\mu_t)$ is an abstract Wiener space, i.e.,
$\mu_t$ is exactly the Gaussian measure from Section 2.
Therefore, if $\dim H=\infty$, by Corollary \ref{cor.2.4}
$$
\mu_t(H+x)=0  \mbox{ for all } x\in E,
$$
hence for all $t>0$
$$
\nu_t(H)=\int\int 1_H(tb+z+y)\, \mu_t(dy)\, \nu_t^{0}(dz)=0.
$$
So, $(H)(ii)$ holds in this case.
\end{exm}

Let $\alpha_n \in (0,\infty)$, $n\in \N$, such that $\alpha_n\nearrow \infty$ as $n\to \infty$ and
$$
\leqno{(3.7)}\quad\quad\quad \quad
\sum_{n=1}^\infty \alpha_n \lambda_n < \infty.
$$
Let us fix $x\in E \setminus H$, and   $e_n^x$, $n \in \N$, be
as in Lemma \ref{lem:3.2}.
Define $q_x: E \to \overline{\R}_+$ by
$$
\leqno{(3.8)}\quad\quad
q_x(z) := \left[ \displaystyle\sum_{n=1}^\infty \alpha_n \lambda_n\,
{}_{E'} \langle e_n , z \rangle _E^2
+ \left( \sum_{n=1}^\infty  2^{-\frac{n}{2}}   | {}_{E'}\langle e^x_n,  z \rangle_{E} |\right)^2
\right]^{\frac{1}{2}},
$$
where $\{ e_n: n\in \N \}$ is the special orthonormal basis of $H$ from above.
Then clearly $q_x$ has compact level sets in $E$. Define again
$$
E_x:= \{z \in E : \,  q_x ( z ) < \infty \} .
$$
Then obviously  $x\not\in E_x$. Furthermore, we have an analog  of
Proposition \ref{prop:p6}.

\begin{prop} \label{prop.3.3}  
Let $t>0$. Then the following assertions hold.

\noindent $(i)$  $q_x\in L^2(E, \nu_t)$, in particular
$\nu_t(E_x)=1$ and $\nu_t(H+x)=0$.

\noindent $(ii)$ $H\subset E_x$ continuously.

\noindent
$(iii)$ For all $z\in E$ we have
$$
||z||\leq q_x (z).
$$
In particular, $(E_x, q_x)$ is complete. Furthermore, $(E_x, q_x)$
is compactly embedded into $(E, ||\cdot  ||)$.
\end{prop}

\begin{pf}
$(i)$  By $(H) (i)$ we have
$$
\leqno{(3.9)}\quad\quad
\displaystyle\int_E q_x^2(z) \nu_t(dz)\leq
C(1+t^2) \sum_{n=1}^\infty \alpha_n \lambda_n+
\left(\sum_{n=1}^\infty 2^{-\frac{n}{2}}
\sqrt{ \int_E  {}_{E'}\langle e^x_n,  z \rangle_{E}^2 \, \nu_t(dz) }\, \right)^2\leq
$$
$$
C(1+t^2)
\left( \sum_{n=1}^\infty \alpha_n \lambda_n +
(\sum_{n=1}^\infty 2^{-\frac{n}{2}})^2 \right) <\infty.
$$

(ii) This is obvious by $(2.1)$ and $(3.7)$.

$(iii)$  By $(3.5)$ we have for all $z \in E$
$$
\leqno{(3.10)}\quad\quad
 \displaystyle\sum_{n=1}^\infty \alpha_n \lambda_n\,
{}_{E'} \langle e_n , z \rangle _E^2 =
 \displaystyle\sum_{n=1}^\infty \alpha_n \lambda^{-1}_n
\langle e_n , z \rangle _E^2 =
 \displaystyle\sum_{n=1}^\infty \alpha_n \langle\overline e_n , z \rangle _E^2 .
$$
Hence since $\alpha_n\nearrow \infty$ as $n\to \infty$, we have
$$
q_x^2(z)\geq \alpha_1 ||z||_E^2,
$$
and, therefore, $(E_x, q_x)$ is complete by Fatou's Lemma and
$(E_x, q_x)$  is compactly embedded  into $(E, ||\cdot  ||)$.

\hfill $\square$
\end{pf}

The following result is an analog to Theorem \ref{thm:t8}
for infinite dimensional L\'evy processes.

\begin{thm}\label{th3.5} 
Assume that $(HS)$ and $(H)$ hold. Let $v_0^x:=U_1 q_x^2$ and for
every $z\in E$, $v_z^x:=v_0^x\circ T_z^{-1}$. Then $v_z^x$ is a
compact Lyapunov function such that  $E_x+z=[v_z^x< \infty]$ and
each $E_x+z$ is invariant with respect to $(P_t)_{t\geq 0}$. In
particular, $E_x+z$ is left invariant by the infinite dimensional
L\'evy process $X=(\Omega,\cf,\cf_t,X_t,\theta_t,P^{x})$.
Furthermore, the restriction of $X$ to $E_x+z$ is c{\`a}dl{\`a}g
in the trace topology.
\end{thm}

\begin{pf}
For $y\in E$, using the sublinearity of $q_x$, by $(3.9)$
we obtain that for some constant $\widetilde C >0$
$$
P_t q_x ^2 (y) \leq 2 q_x ^2 (y) + 2 \int _E q_x ^2 ( z ) \nu _t ( dz )  \leq
2q_x^2(y)+2 \widetilde C ( 1 + t ^2 ).
$$
Hence
$$
\leqno{(3.11)} \quad\quad
v_0^x (y) = U_1 q_x^2 (y) = \displaystyle\int _0 ^\infty e ^{-t} P_ t q_x^2 (y) dt \leq 2 q_x ^2 (y)
+ 2 \widetilde C \int _0 ^ \infty (1 + t ^2 ) e ^{-t} dt.
$$
On the other hand, since $q_x$ is a norm, for all $y,z\in E$
by the triangle inequality we have that
$$
q_x^2(y+z)\geq ( q_x(y) - q_x(z) )^2\geq \frac{1}{2} q_x^2(y) - q_x^2(z).
$$
Hence by $(3.9)$
$$
P_t q_x^2(y)\geq \frac{1}{2} q_x^2(y)- \int_E q_x^2(z)\nu_t(dz)\geq
 \frac{1}{2} q_x^2(y)- \widetilde{C} (1+t^2)
$$
and therefore
$$
\leqno{(3.12)} \quad\quad \displaystyle
v_0^x(y)=U_1 q_x^2(y)=\int_0^\infty e^{-t} P_t \, q_x^2(y)\, dt\geq
 \frac{1}{2} q_x^2(y)- \widetilde{C} \int_0^\infty e^{-t} (1+t^2)\, dt.
$$
Finally, by $(3.11)$  and $(3.12)$ it follows that
$$
E_x = [v_0 < \infty ].
$$
$v_0^x$ is a Lyapunov function for $( P_t ) _{t \geq 0}$, which is
compact by $(3.12)$.

\noindent Since the measure $\nu_t$ is carried by $E_x$, it
follows by the same argument as in the proof of  Lemma
\ref{lem:l7} that each $E_x+z $ is an invariant set for
$(P_t)_{t\geq 0}$.

To prove the next part of the assertions let us more generally
consider any set $L\in\cb(E)$ instead of $E_x+z$ just with the
property that $P_t 1_L=1_L$ for all $t>0$. Then $1_L\in\ce(\cu)$,
hence it is finely continuous and therefore
$[1_L=0]=[1_L<\frac{1}{2}]$ is finely closed and finely open.
Consequently, for all $x\in L$, $t>0$,
$$
P^x([1_L(X_t)>\frac{1}{2}])=P^x([X_t\in L ]
)=E^x{[1_L(X_t)]}=P_t1_L(x)=1
$$
and thus, since $t\longmapsto 1_L(X_t)$ is continuous because
$1_L\in\ce(\cu)$, we obtain
$$
P^x({X_t\in L} \ \ \forall t \geq 0)=P^x(\bigcap_{t \geq
0}[1_L(X_t)>\frac{1}{2}]) =P^x(\bigcap_{t\in \Q^{+}} [
1_L(X_t)>\frac{1}{2}])=1.
$$
To prove the final assertion  let $X'$ be the restriction of $X$
to $L$,  $\cu'=(U'_\alpha)_{\alpha>0}$ be its resolvent, and
recall that $U_\alpha(C_b(E))\subset C_b(E)$ for all $\alpha >0$,
where $\cu=(U_\alpha)_{\alpha>0}$ is the resolvent of $X$.
Consequently, $U'_\alpha$    maps $C_b(E)|_L$ into $C_b(E)|_L$ for
all $\alpha>0$. From the first part of the proof there exists on
$L$ a real valued compact Lyapunov function with respect to
$\cu'$. The claimed c{\`a}dl{\`a}g property of $X'$ follows now by
Theorem 5.2 from \cite{BeRo09}.
\hfill $\square$
\end{pf}

\begin{rem}     
$(i)$ The analog of Remark \ref{rem2.7'} holds, i.e. $E$ is an
uncountable disjoint union of Borel sets which are invariant for
the L\'evy process on $E$.

$(ii)$ Subsection $3.2$ from \cite{BeRo11} presents an informal
description of constructing compact Lyapunov functions for the
infinite dimensional L\'evy processes.
\end{rem}

\begin{exm}\label{exm3.6} 
   Let $(S,\mcB, \sigma)$ be a finite measure space and
   $H:= L ^2 ( S , \mcB, \sigma)$.
   Define $\l :H \to \C$ by
   \[\l (h ) := \int _S ( 1- e^{ih}) d \sigma , \;  h \in H.\]
   Then  $\l (0) =0$, $\l$ is negative definite and continuous on $H$.
   Choosing a Hilbert-Schmidt extension
   $E$ of $H$ as above there exist probability measures
   $\nu_t, \; t > 0,$ on $(E, \mcB(E))$ such that
$$
\widehat \nu_t ( \xi ) = \int_E e^{i {}_{E'}\langle \xi, z \rangle_E } \,  \nu_t ( dz)=
e^{- t \int _S ( 1 - e^{i\xi }) d\sigma}, \; \xi \in E'.
$$
   $\nu_t $ is just the Poisson measure with intensity $t$ on $E$.
   Hence for all $\xi \in E'$
   \[\int \bracket{\xi, z}^2 \nu_t ( d z) =
   t \int _S \xi ^2 d \sigma + t ^2 \left(\int _S \xi d \sigma \right)^2
   \leq \sup(2 \sigma (S)^2) ( 1 + t^2 ) | \xi |_H ^2. \]
   In particular, (H)(i) holds.
   \end{exm}

Now take $S=(0,1)$, $\cb$=Borel $\sigma$-algebra on $(0,1)$ and
$\sigma$=Lebesgue measure $ds$.
Let $H_0^1$ be the Sobolev space of order $1$ in $L^2((0,1),ds)$
with Dirichlet boundary conditions.
Let
$$
E:=(H_0^1)'(=H^{-1}).
$$
Then we have the Hilbert-Schmidt embeddings
$$E'=H_0^1\subset L^2((0,1),ds):=H\subset E.
$$
So, each $\nu_t$ extends to a probability measure on $(E,\cb(E))$.
Since $H_0^1$ continuously embeds into the bounded continuous on
$(0,1)$ equipped with the sup-norm, $E$ contains all measure of
finite total variation. It is, however, well-known (see, e.g.,
\cite{Kun99}) that each $\nu_t$ is supported by positive measures
of type $\ds\sum_{i=1}^{N}\varepsilon_{x_i}$, where
$\varepsilon_{x_i}$ is a Dirac measure with mass in $x_i\in
[0,1]$, $1\leq i\leq N_{x}\in \N$, and $x_i$ are pairwise
distinct. In particular, $\nu_t(H)=0$ for all $t>0$.
So, also $(H) (ii)$ holds in this case.

Similar arguments can be used in the case where $S$ is replaced by
an open bounded set in $\R^d$. Then one has to take $E$ as the
dual of a Sobolev space of sufficiently (with respect to $d$) high
order. Likewise one can treat the case $S=\R^d$, but then one has to
use weighted Sobolev spaces.

\section{Potential theory} 

\subsection{Preliminaries}

In this section we consider the Banach space $E$
and the Hilbert space $H$ as in Section 2.
Let $(\nu_t)_{t \geq 0}$ be a convolution semigroup
of probability measures
on $(E,\cb)$ and $(P_t)_{t \geq 0}$ the associated
family of Markovian kernels:
$$
P_t f(x)=\int_{E}f(y)p_t(x,dy)=\int_{E}f(x+y)\nu_t(dy),\;\;\;f\in
p\cb,\;x\in E,
$$
where $p_t(x,\cdot)$ is the probability measure on
$(E, \cb)$ such that
$$
p_t(x,A):=\nu_t(A-x) \mbox{ for all }A\in\cb .
$$

Let further $\cu=(U_{\alpha})_{\alpha>0}$ be the Markovian
resolvent of kernels on $(E,\cb)$ associated with $(P_t)_{t \geq
0}$, i.e., $\ds U_{\alpha}:=\int_{0}^{\infty}e^{-\alpha t}P_t \
dt,$ $\alpha>0$, and set $U:=\int_{0}^{\infty}P_t dt$. $U$ is
called potential kernel of $\cu$. Clearly, for
$\cu_{\beta}:=(U_{\beta+\alpha})_{\alpha >0}$ the corresponding
potential kernel is $U_{\beta}$.

We consider an orthonormal basis $\{e_n:n\in \N^{*}\}$ of $H$
formed by $e_n\in E{'}$, $n \in \N^{*}$. For each $n$ define
$$
\widetilde{P_n}:E\rightarrow
H_n:=span\{e_1,e_2,\cdots,e_n\}\subset E{'} \subset H
$$
by
$$
\widetilde{P_n}z:=\sum_{k=1}^{n} {}_{E'}\langle e_k,z \rangle_{E} \ e_k,\;\;\;z\in E.
$$

Whenever necessary, $H_n$ is identified with $\R^{n}$.
For each $t>0$ and $n\in\N^{*}$
we consider the probability measure $\nu_{t}^{\{n\}}$
on $\R^{n}$ defined by
$$\nu_{t}^{\{n\}}:=\nu_t\circ\widetilde{P_n}^{-1}.
$$
Analogously, we consider the kernel $P_{t}^{\{n\}}$ on
$(\R^{n},\cb(\R^{n}))$ induced by $\nu_{t}^{\{n\}}$:
$$
P_{t}^{\{n\}}\varphi(x)=
\int_{\R^n}\varphi (x+z)\nu_{t}^{\{n\}}(dz),\;\;\;\varphi\in p\cb(\R^n),\;x\in\R^{n}.
$$
We obtain a Markovian semigroup of kernels
$(P_t^{\{n\}})_{t \geq 0}$ on $(\R^{n},\cb(\R^{n}))$ and let
$\cu^{n}=(U_{\alpha}^{\{n\}})_{\alpha >0}$ be the associated resolvent of kernels.

Let $n\in\N^{*}$, $t>0$, and $f$ be a positive cylinder function on $E$
based on $H_n$, i.e.,
there exists a function $\varphi\in p\cb(\R^n)$
such that   $f=\varphi\circ\widetilde{P_n}$.
Then for all $x\in E$ we have
$$
P_t f(x)=\int_{E}f(x+y)\nu_t(dy)=
\int_{\R^n}\varphi(\widetilde{P_n}x+z)\nu_t^{\{n\}}(dz)=P_t^{\{n\}}\varphi(\widetilde{P_n}x).
$$
Consequently,  for all $\alpha>0$ we have

$$ \leqno{(4.1)} \hspace*{34mm} U_{\alpha}f =
(U_{\alpha}^{\{n\}}\varphi)\circ \widetilde{P_n}.
$$\\

\begin{prop} \label{prop4.1}
Let $v \in p \cb(\R^n)$ and $\beta>0$.
Then $v$ is $\cu_{\beta}^{\{n\}}$-excessive
(resp. $\cu_{\beta}^{\{n\}}$-supermedian, i.e.,
$\alpha U_{\beta+\alpha}^{\{n\}}v \leq v$ for all $\alpha>0$) if and only if
$v\circ\widetilde{P_n}$ is $\cu_{\beta}$-excessive
(resp.  $v\circ\widetilde{P_n}$ is $\cu_{\beta}$-supermedian).
\end{prop}

\begin{pf}
The assertion follows from the equality $(4.1)$:
$$
U_{\alpha}(v \circ \widetilde{P_n})=(U_{\alpha}^{\{n\}}v)\circ \widetilde{P_n}.
$$
\hfill $\square$
\end{pf}

We assume further that $(P_t)_{t \geq 0}$ (resp. $(P_t^{\{n\}})_{t
\geq 0}$) is the transition function of a right process
$X=(\Omega,\cf,\cf_t,X_t,\theta_t,P^{x})$ with state space $E$
(resp.
$X^{\{n\}}=(\Omega^{\{n\}},\cf^{\{n\}},\cf_t^{\{n\}},X_t^{\{n\}},$ $\theta_t^{\{n\}},P^{x})$
with state space $\R^n$), i.e.,
$$
P_t f(x)=E_{x}(f\circ X_t),\;\;\;x\in E,\; f\in p\cb(E).
$$

\begin{rem} \label{rem.4.2}  
$(i)$  The Gaussian measures in an abstract
Wiener space (presented in Section 2) and the convolution
semigroup of a L\'evy process on a Hilbert space (studied in
Section 3) are examples for which the results from this section
apply.

$(ii)$ If $\nu_t=\mu_t$, a Gaussian measure with parameter $t$ in
an abstract Wiener space, then $\nu_t^{\{n\}}$ is the
$n$-dimensional Gaussian measure with parameter $t$. Consequently,
Proposition 4.1 has the following interpretation: {\rm every
superharmonic function in an $n$-dimensional Euclidean space is
"superharmonic" with respect to the Gross-Laplace operator, i.e.,
it is an excessive function for the infinite dimensional Brownian
motion, when it is canonically transported on the abstract Wiener
space.}
\end{rem}

\begin{cor}     
Suppose that $(\nu_t)_{t \geq 0}$ is the convolution semigroup of
a L\'evy process on an Hilbert space as in Section 3. If for some
$n\in \N^*$ the process $X^{\{ n\} }$ is transient then $X$ is
also transient. If $X$ is not transient then $X^{\{ n\} }$ is
recurrent for all $n$.
\end{cor}

\begin{pf}
If the process $X^{\{ n\} }$ is transient, or equivalently  the
potential kernel $U^{\{ n\} } = \int_0^\infty  P^{\{ n\} }_t dt $
of $X^{\{ n\} }$ is proper, then by $(4.1)$ we get that the
potential kernel $U$ of $X$ is also proper. The second assertion
follows from the first one and by the transience--recurrence
dichotomy which holds for L\'evy processes (cf., e.g., Theorem
35.4 in \cite{Sa99}).
\hfill $\square$
\end{pf}

\subsection{Excessive measures and the energy functional} 

Let $Exc(\cu)$ be the set of all {\it $\cu$-excessive
measures} on $E$: $\xi\in Exc(\cu)$ if and only if it is a
$\sigma$-finite measure on $(E,\cb)$ such that $\xi\circ \alpha
U_{\alpha}\leq \xi$ for all $\alpha>0$.

By $Pot(\cu)$ we denote the set of all {\it potential}
$\cu$-excessive measures, i.e. all $\sigma$-finite measures $\xi$ of the form
$\xi=\mu\circ U$, where $\mu$ is a  measure on
$(E,\cb)$. Clearly, by the resolvent equation we have that
$Pot(\cu)\subset Exc(\cu)$. Note that the {\bf mass uniqueness
principle} holds for the Gaussian measures in an abstract Wiener
space  and the convolution semigroup of a L\'evy process on a
Hilbert space:

{\it If $\beta>0$ and $\mu, \nu$ are two positive measures on
$(E,\cb)$ such that $\mu\circ U_\beta$, $\nu \circ U_\beta$ are
$\sigma$-finite and  $\mu\circ U_\beta= \nu \circ U_\beta$, then
$\mu=\nu$.}

The assertion follows from $(10.40)$ in \cite{Sh88}; see
Proposition 5 in \cite{Ca80} for the Gaussian case.\\

If $\beta>0$ then the {\it energy functional}
$L_{\beta}:Exc(\cu_{\beta})\times \ce(\cu_{\beta})\longrightarrow
\overline{\R}_{+}$ is defined by
$$
L_{\beta}(\xi,v):=\sup\{\mu(v):Pot(\cu_{\beta})\ni\mu\circ
U_{\beta}\leq \xi\}.
$$

The following result is  a consequence of $(4.1)$ and Proposition
\ref{prop4.1}.

\begin{cor}      
The following assertions hold.

$(i)$ If $\xi\in Exc(\cu_{\beta})$ then
$\xi\circ\widetilde{P_n}^{-1}\in Exc(\cu_{\beta}^{\{n\}})$
provided it is a $\sigma$-finite measure on $\R^{n}$. If in
addition $\xi\in Pot(\cu_{\beta})$ then
$\xi\circ\widetilde{P_n}^{-1}\in Pot(\cu_{\beta}^{\{n\}})$.

$(ii)$ Let $\xi\in Exc(\cu_{\beta})$ such that
$\xi\circ\widetilde{P_n}^{-1}$ is  $\sigma$-finite, $v\in
\ce(\cu_{\beta}^{\{n\}})$, and let $L_{\beta}^{\{n\}}$ be the
energy functional with respect to $\cu_{\beta}^{\{n\}}$. Then
$$
L_{\beta}^{\{n\}}(\xi\circ\widetilde{P_n}^{-1},v)=L_{\beta}(\xi,v\circ\widetilde{P_n}).
$$
\end{cor}

\subsection{Absence of a reference measure} 

Recall that a right Markov process satisfies the {\it hypothesis
(L)} of P.A. Meyer provided that there exists a finite measure on
$(E, \cb )$ with respect to which  all the measures $U_\alpha (x,
\cdot),$  $x\in E,$ are absolutely continuous, where
$\cu=(U_\alpha)_{\alpha>0}$ is the resolvent family of the
process. Such a measure is called {\it reference measure} for
$\cu$.  Recall that the {\it fine topology} is the topology on $E$
generated by $\ce(\cu_\beta)$.

\begin{prop}   
The hypothesis $(L)$ of P.A. Meyer does not hold for the L\'evy
processes on an infinite dimensional Hilbert space.
\end{prop}

\begin{pf}
The main argument in the proof is the same as in the Gaussian case
(cf. Proposition 8 in \cite{Ca80}), namely, the existence of an
uncountable family of mutually disjoint finely open sets. More
precisely, assume that there exists a reference measure $\lambda$
for $\cu$. Note that $\lambda$ charges every non-empty finely open
set. Indeed, if $G\in \cb$ is finely open and we suppose that
$\lambda (G)=0$ then $U_\beta(1_G)\equiv 0$, which contradicts the
fact that $U_\beta(1_G)(x)>0$ for all $x\in G$. (cf., e.g.,
Proposition 1.3.2 from \cite{BeBo04}). Since $\dim H=\infty$,
there exists $x\in E\setminus H$ and the space $E_x$ defined in
Section 3.  By Theorem \ref{th3.5} the sets $E_x+z$, $z\in E$, are
invariant with respect to $(P_t)_{t\geq 0}$. In particular,
$E_x+z$ is finely open for every $z\in E$. Because  $x\not\in
E_x$, it follows that $(E_x + {\alpha x})_{\alpha \in \R_+}$ is an
uncountable  family of mutually disjoint sets and from the above
considerations we get $\lambda (E_x+\alpha x)>0$ for all $\alpha
\in \R_+$, which leads to a contradiction.
\hfill $\square$
\end{pf}

\subsection{Reduced functions and polar sets}  

If $M \subset E$ and $v
\in\ce(\cu_{\beta})$, then the {\it reduced function} (with
respect to $\cu_{\beta}$) of $v$ on $M$ is
the function $R_{\beta}^{M}v$ defined by:

$$
R_{\beta}^{M}v:=\inf\{u\in\ce(\cu_{\beta}):u\geq v \mbox{ on }M\}. 
$$

If $M$ is a Souslin subset of $E$ then the reduced function
$R_{\beta}^{M}v$ is universally $\cb$-measurable.
The maps $v\longmapsto R^M_\beta v$ and  $v\longmapsto \widehat{R^M_\beta} v$ 
extend to kernels on $E$ and by Hunt's Theorem we have
$$
R^{M}_{\beta}v(x)=E^{x}(e^{-\beta D_M} v \circ X_{D_M};D_M<\infty),
$$
$$
\widehat{R_{\beta}^{M}}v(x)=E^{x}(e^{-\beta T_M}v\circ X_{T_M};T_M<\infty),
$$
where
 $D_M(\omega):=\inf \{t \geq 0\;|\;X_t(\omega)\in M\}$,
 $T_M(\omega):=\inf\{t>0\;|\;X_t(\omega)\in M\}$, $\omega\in\Omega$,
and for a $\cu_{\beta}$-supermedian function $u$, $\widehat{u}$
denotes its $\cu_{\beta}$-excessive regularization,
$\ds\widehat{u}(x)=\sup_{\alpha>0}\alpha U_{\beta+\alpha}u(x)$ for
all $x\in E$.

The set $M\in\cb$ is called {\it polar} (resp. {\it $\nu$-polar};
where $\nu$ is a $\sigma$-finite measure on $(E,\cb)$) if
$\widehat{R_{\beta}^{M}}1=0$ (resp. $\widehat{R_{\beta}^{M}}1=0$
$\nu$-a.e.). By the above mentioned Hunt's Theorem a set $M\in \cb$
will be polar (resp.  $\nu$-polar) if and only if
$T_M=\infty$ $P^x$-a.s. for all $x\in E$
 (resp. $T_M=\infty$ $P^\nu$-a.e.).

\begin{cor} \label{cor.4.5}  
If $M\in\cb$, $n\in \N^*$, and $v\in\ce(\cu_{\beta}^{\{n\}})$ then
$$
R_{\beta}^{M}(v \circ \widetilde{P_n})\leq
({}^{\{n\}}\!\! R_{\beta}^{\widetilde{P_n}(M)}v)\circ\widetilde{P_n},
$$
where for a set $F \subset \R^{n}$ we have denoted by
${}^{\{n\}}\!\! R_{\beta}^{F}v$ the reduced  function
(with respect to $\cu_{\beta}^{\{n\}}$) of $v$ on $F$.
In particular, if $\widetilde{P_n}(M)$ is a polar subset
of $\R^n$ then $M$ is a polar subset of $E$.
\end{cor}

\begin{pf}
Let $u\in\ce(\cu_{\beta}^{\{n\}})$, $u \geq v$ on $\widetilde{P_n}(M)$.
Then $u\circ \widetilde{P_n} \geq v\circ \widetilde {P_n}$ on $M$ and
by Proposition 4.1 we have $u \circ \widetilde{P_n} \in \ce(\cu_{\beta})$.
Consequently, we get that
$u\circ \widetilde{P_n} \geq R_{\beta}^{M}(v \circ \widetilde{P_n})$
on $E$ and thus for all $x\in E$ we have
\begin{center}
${}^{\{n\}}\!\!R_{\beta}^{\widetilde{P_n}(M)}v(\widetilde{P_n}x)=
\inf\{u(\widetilde{P_n}x)\, : \, u\in\ce(\cu_{\beta}^{\{n\}}),\;u
\geq v$ on $\widetilde{P_n}(M)\}$ $\geq R^{M}_{\beta}(v \circ
\widetilde{P_n})(x).$
\end{center}

Assume now that $\widetilde{P_n}(M)$ is a polar subset of $\R^{n}$.
Using $(4.1)$ we get for all $x\in E$
$$
U_{\alpha}^{\{n\}}({}^{\{n\}}\!\! R_\beta ^{\widetilde{P_n}(M)}v)(\widetilde{P_n}x)=
U_{\alpha}({}^{\{n\}}\!\! R_\beta ^{\widetilde{P_n}(M)}v\circ\widetilde{P_n})(x)\geq
U_{\alpha}(R^{M}_{\beta}(v \circ \widetilde{P_n}))(x)
$$
and therefore, taking $v=1$ we have
$$
0=\widehat{{}^{\{n\}}\!\! R^{\widetilde{P_n}(M)}}1(\widetilde{P_n}x)\geq\widehat{R_{\beta}^{M}}1(x),
$$
hence $M$ is a polar subset of $E$.
\hfill $\square$
\end{pf}

\begin{prop}\label{prop.4.6}  
Assume that $(\nu_t)_{t \geq 0}$ is the convolution semigroup
of a L\'evy process on an Hilbert space as in Section 3 and
suppose that for all $t>0$ $\nu_t$ charges no proper closed linear subspace of $E$.
Then the points of $E$ are polar sets.
\end{prop}

\begin{pf}
By Corollary \ref{cor.4.5}  it is sufficient to show that the points are
polar for one finite dimensional projection
$(\nu_t^{\{n\}})_{t\geq 0}$ of $(\nu_t)_{t \geq 0}$.
By Theorem 4 in \cite{Br71} it follows that the points are polar for
a L\'evy process in $\R^{n}$, $n \geq 2$, provided that the points
are not finely open sets for all $1$-dimensional projections.
Suppose that $\{0\}\subset \R$ is a finely open set for
$(\nu_t^{\{1\}})_{t \geq 0}$.
Proposition 4.1 implies that $\widetilde{P_n}^{-1}(G)$ is
a finely open subset of $E$ for every $G \subset \R^{n}$
which is finely open with respect to $\cu_{\beta}^{\{n\}}$.
Consequently, the set $F:=\widetilde{P_1}^{-1}(\{0\})$
will be a closed proper subspace of $E$ which is finely open,
hence $U_{\beta}(1_F)>0$ on $F$.
This contradicts the hypothesis on $\nu_t$ which implies $\nu_t(F)=0$.
Therefore $\{0\}\subset \R$ is not finely open and we conclude
that the set $\{0\}\subseteq E$ is polar.
\hfill $\square$
\end{pf}

\begin{prop}\label{prop.4.7} 
Let  $(\nu_t)_{t \geq 0}$ be either the Gaussian semigroup in an
abstract Wiener space or the convolution semigroup of a L\'evy
process on an Hilbert space as in Section 3, satisfying hypotheses
$(HS)$ and $(H)$. Then the  "Cameron-Martin" space $H$ is a polar
set.
\end{prop}

\begin{pf}
Let $x\in E \setminus H$. By Corollary \ref{cor.2.4} and Lemma
\ref{lem:l7}  (in the Gaussian case) and by Proposition
\ref{prop.3.3} (in the L\'evy process case) there exists
$E_x\in\cb$, a linear subspace of $E$, such that $H\subset E_x$,
$\nu_t(E_x)=1$ and $x\notin E_x$. Using again Lemma \ref{lem:l7}
(in the Gaussian case) and Theorem \ref{th3.5} in the L\'evy
process case) we get  that $E_x $ is invariant with respect to
$(P_t)_{t\geq 0}$, hence $1_{E_x}\in\ce(\cu_{\beta})$.
Consequently, we get $R_{\beta}^{H}1(x)\leq 1_{E_x}(x)=0$ and thus
$R_{\beta}^{H}1=0$ on $E\setminus H$. Since $p_t (y, H)=0$ for all
$y\in E$ and $t>0$, we get $U_\alpha(1_H)=0$ and so
$$
\widehat{R_{\beta}^{H}}1(x)=
\lim_{\alpha\rightarrow\infty}\alpha U_{\alpha}(R^{H}_{\beta}1)(x)=0 \mbox{ for all } x\in E .
$$
\hfill $\square$
\end{pf}

\begin{rem} 
$(i)$ The result of Proposition \ref{prop.4.7} was proved in the
Gaussian case in \cite{Ca80}, Proposition 4. Note that the main
probabilistic argument used in that proof (see Remark 7 in
\cite{Ca80}) remains valid here: The property of $E_x+x$ to be
invariant with respect to $(P_t)_{t\geq 0}$ implies that the
process starting from $x$ never leaves the set $E_x+x$. Since
$H\subset E\setminus (E_x+x)$, it follows that the process
starting from $x$ never hits $H$.

$(ii)$ If $H$ is polar, then clearly all the points are polar sets.
So, the conclusion of Proposition \ref{prop.4.7}
is stronger than that
of Proposition \ref{prop.4.6}.
\end{rem}

\subsection{Choquet capacities and quasi continuity} 

In this subsection we assume
again that $(\nu_t)_{t \geq 0}$ is the convolution semigroup of a
L\'evy process on an Hilbert space as in Section 3;
see Theorem \ref{thm:2.7} and  \cite{BeCoRo08}  for the Gaussian case.

In Remark \ref{rem.2.9} $(ii)$ we recalled Carmona's question on the
existence of a relevant capacity for the infinite dimensional Brownian motion.
We can present now the corresponding capacity for the L\'evy processes.
Note that in this case, since these processes are not necessarily   transient,
we have to consider the "$\beta$-level" capacity, $\beta>0$.

Let $p:= U_\beta f_0$, with $0<f_0\leq 1$, $f_0\in p\mathcal B$,
and let $\lambda$ be a finite measure on $(E,\mathcal B)$. Then
the functional $M\longmapsto c_\lambda (M)$, $M\subset E$, defined by
$$
c_\lambda(M) := \inf \bigl\{ \lambda(R_\beta^G p) \, : \,
M\subset G\;\text{open} \bigr\}
$$
is a Choquet capacity on $E$ (see e.g.\ \cite{BeBo04}).

We complete this subsection with an  analog of Theorem
\ref{thm:2.7} for L\'evy processes.

\begin{thm} 
$(i)$ The topology of $E$ is a Ray one and the capacity
$c_\lambda$ is {\it tight}, i.e., there exists an increasing
sequence $(K_n)_n$ of compact sets such that $\inf_n c_\lambda
(E\setminus K_n)=0$.

$(ii)$  Let $M\in\mathcal B$. Then
$$
c_\lambda (M) = \lambda(R_\beta^M p)  = \sup \bigl\{ \nu(p\cdot
1_M) \, :\,  \nu\circ U_\beta \leq \lambda \circ U_\beta
\bigr\}\;.
$$
The set $M$ will be $\lambda$-polar and $\lambda$-zero  if and
only if $c_\lambda(M)=0$.

$(iii)$  Every $\cu_\beta$-excessive function of the form $U_\beta
f$, $f\in p\cb$,  is $c_\lambda$-quasi continuous, provided it is
finite $\lambda$-a.e. More generally, every $(\beta)$-level
potential of a continuous additive functional (cf. \cite{Sh88} or
\cite{BeBo04} in the transient case) is $c_\lambda$-quasi
continuous if it is finite $\lambda$-a.e. In particular, every
$\cu_\beta$-excessive function is $c_\lambda$-quasi lower
semicontinuous.
\end{thm}

\begin{pf}
$(i)$ Let $C_{bl}(E)$ be the set of all bounded Lipschitz
continuous functions on $E$. Using $(3.4)$ one can check that
$(U_\alpha)_{\alpha>0}$ induces a strongly continuous resolvent of
contractions on $C_{bl}(E)$ and then one can construct an
appropriate Ray cone (see Proposition $2.2$ from \cite{BeRo09} for
details). The tightness property follows by \cite{MaRo92} (see
also \cite{BeBo05}) since we already remarked in Section $3$ that
an infinite dimensional L\'evy process has c\`adl\`ag paths.

Assertion $(ii)$ is a consequence of  Proposition 1.6.3 and
Proposition 1.6.4 from \cite{BeBo04}, because by $(i)$ the
topology of $E$ is a Ray one.

As in the proof of Theorem \ref{thm:2.7}, assertion $(iii)$
follows by  Proposition $3.2.6$ from \cite{BeBo04}, using again
the property of the topology to be a Ray one.
\hfill $\square$
\end{pf}

\subsection{Existence of bounded invariant functions}

\begin{rem}  
$(i)$ Suppose that $(\nu_t)_{t \geq 0}$ is the convolution
semigroup of a L\'evy process on an infinite dimensional Hilbert
space as in Section 3 and $x\not\in H$. By Theorem \ref{th3.5} the
function $1_{E_x}$ is invariant with respect to $(P_t)_{t\geq 0}$,
it is identically equal to one on $H$ and zero at $x$. This  shows
that the answer given by R. Carmona (see Remark 6 in \cite{Ca80})
to a conjecture of V. Goodman (cf. \cite{God73}, page 219) for the
infinite dimensional Brownian motion, remains valid for the L\'evy
processes on an Hilbert space.

$(ii)$ Unbounded invariant functions may be further constructed as in  \cite{Ca80},
the proof of Proposition 3, namely, consider the function $f$ defined as
$$
f= \sum_{n=1}^\infty r_n 1_{\frac{1}{n}x+E_x},
$$
where $(r_n)_n$ is a sequence of real numbers with $\lim_{n\to \infty} r_n=\infty$.
Then  clearly $f$ is invariant and it is unbounded in every neighborhood of each point.

$(iii)$  Let $v\in bp\cb$ be invariant with respect to $(P_t)_{t\geq
0}$, assume that the L\'evy process has continuous paths
(i.e., $M$ in $(3.3)$ is the zero measure),
and consider an open set $V\subset E$ which is transient, i.e., we
have a.s. $\sup \{ t>0 \, : \, X_t\in V \}< \infty$. Then the
function $v$ is harmonic on $V$ in the sense considered in the
Gaussian case (see Section 5 below): $v$ is finely continuous and
there exists $\rho>0$ such that
$$
v(x) = P_{T_{E\setminus B_r(x)}} v(x)
$$
for all $r<\rho$ whenever $\bar B_r(x) \subset V$; $\bar B_r(x)$
denotes the closed ball or radius $r$ centered at $x$. Indeed,
since $V$ is transient we get that a.s. ${T_{E\setminus
B_r(x)}}<\infty$. The assertion follows from a straightforward
consequence of Dynkin's formula (cf., e.g.,  $(12.18)$  in
\cite{Sh88}): if $v$ is a bounded $\cu$-invariant function and $T$
is a terminal time with $T<\infty$ a.s., then $v=P_Tv$.
\end{rem}

\subsection{Domination principle} 

\begin{prop} \label{prop:domin}
Let $\mu$, $\nu$ be two $\sigma$-finite measures on $(E, \cb)$,
  $G\in \cb$  a finely open set such that
$\mu(E\setminus G)=0$. Assume that $\mu\circ U_\beta$, $\nu\circ
U_\beta$ are $\sigma$-finite  measures and $\mu\circ U_\beta \leq
\nu\circ U_\beta$ on $G$ for some $\beta>0$. Then $\mu\circ
U_\beta \leq \nu\circ U_\beta$ on $E$.
\end{prop}

\begin{pf}
For $\xi\in Exc(\cu_\beta)$ and $M\in \cb$ define ${}^*\! R^M \xi
:= \bigwedge \{ \eta \in Exc(\cu_\beta)\, : \, \eta\geq \xi$  on
$M \}$, where $\bigwedge$ denotes the infimum in $Exc(\cu_\beta)$.
If $u\in \ce(\cu_\beta)$, then by Theorem 1.4.12 in  \cite{BeBo04}
$$
\leqno{(4.2)} \quad\quad\quad L_\beta( {}^*\! R^G \xi , u)=
L_\beta (\xi, R^G_\beta u).
$$
Since $R^G_\beta U_\beta f=U_\beta f$ on $G$, $f\in bp\cb$, and
using $(4.2)$ we have
$$
\mu\circ U_\beta(f)= \mu(R^G_\beta U_\beta f)= L_\beta ({}^*\!
R^G_\beta(\mu \circ U_\beta), U_\beta f)= {}^*\! R^G_\beta(\mu
\circ U_\beta)(f) \leq \nu\circ U_\beta(f).
$$
\hfill $\square$
\end{pf}

\begin{rem}
Proposition \ref{prop:domin} is a version of the domination
principle stated for the Gaussian case in Proposition 6 from
\cite{Ca80}. However, our statement is valid for general right
processes, it  holds also for $\beta =0$ in the transient case
(i.e., if the kernel $U= \int_0^\infty P_t \, dt$ is proper), and
it is closer to the original assertion from \cite{Hu57}. The use
of the "duality formula" $(4.2)$ enabled us to avoid the
assumption on the strong duality  from \cite{Hu57}.
\end{rem}

\subsection{Balayage principle}  

The next proposition points out that the balayage principle holds for
the infinite dimensional L\'evy processes; see Proposition 7 in \cite{Ca80}
for the Gaussian case.

\begin{prop}\label{balayageprinc} 
Let $\beta>0$,  $M\in \cb$, $\nu$ a $\sigma$-finite measure on
$(E, \cb)$, and consider the measure $\nu_M$ defined by
$$
\nu_M := \nu \circ \widehat{R^M_\beta} .
$$
Then $\nu_M$ is carried  by the fine closure of $M$,
$\nu_M \circ U_\beta\leq \nu \circ  U_\beta$, and
$$
\nu_M\circ U_\beta = \nu \circ U_\beta \mbox{ on } M.
$$
\end{prop}

\begin{pf}
By Proposition 1.7.11 from \cite{BeBo04}  the measure $\nu_M$ is carried  by the fine closure of $M.$
Since $\widehat{R^M_\beta} u \leq u$ for every $u\in \ce(\cu_\beta)$,
it follows that $\nu_M\circ U_\beta\leq \nu \circ U_\beta$.
If $\cb\ni F\subset M$ then $\widehat{R^M_\beta} U_\beta (1_F)= U_\beta (1_F)$ and so
$\nu_M\circ U_\beta (F)= \int_E \widehat{R^M_\beta} U_\beta (1_F)\, d\nu = \nu \circ U_\beta (F)$.
\hfill $\square$
\end{pf}

\begin{rem} \label{rem.4.11}  
 $(i)$  The assertion of Proposition \ref{balayageprinc}  holds also for
$\beta =0$ in the transient case.

$(ii)$ Recall that the fine closure of $M$ is precisely the union of $M$ with the set of all its regular points;
a point $x\in E$ is called {\rm regular} for $M$ if $P^x(T_M=0)=1$
(see, e.g.,  \cite{Sh88} or \cite{BeBo04}).

$(iii)$ The measure $\nu_M$ is called the {\rm balayage of $\nu$ on
$M$}. Proposition \ref{balayageprinc} offers  an analytic
construction of the balayage of a measure, and therefore,  in the
particular case of the Brownian motion on an abstract Wiener
space,  this gives the answer to a question of  R. Carmona (cf.
Remark 8 in \cite{Ca80}).
\end{rem}

\noindent {\bf Open problem:} It is still open the question
(formulated in \cite{Ca80}, page 38)  whether the axiom of
polarity holds for the infinite dimensional Brownian motion.

\section{Dirichlet problem and controlled convergence}  

Let $\mathcal W = (\Omega,\mathcal F,\mathcal F_t, W_t,\theta_t,P^x)$
be the path continuous Borel right process with state space $E$, having
$(P_t)_{t\geq 0}$ as transition function, given by Theorem \ref{thm:2.7};
recall that $\mathcal W$ is called the \emph{Brownian motion} on $E$.

We already noted in Section 2  that the process $\mathcal
W$ is \emph{transient}, i.e.,
there exists a bounded strictly positive
$\mathcal B$-measurable function $f$ such that $Uf=\int_0^\infty P_t f dt$ is finite.
Therefore in this case we may use the "0-level" excessive functions and potential theoretical tools.
Let $M\in\mathcal B$ and $P_{T_M}$ be the
associated \emph{hitting kernel},
$$
  P_{T_M} f(x)
  = E^x(f\circ W_{T_M}\,;\, T_M<\infty)\;,
  \quad x\in E,\; f\in p{\mathcal B},
$$
where $T_M(\omega) := \inf \bigl\{ t>0 : \,  W_t(\omega)\in M
\bigr\}$, $\omega\in \Omega$. If $u\in\mathcal E(\mathcal U)$,
then $P_{T_M}u=\widehat{R^M}u$.

\begin{rem} \label{rem:hitting} 
If $V$ is an open set and $x\in V$  then the hitting distribution
$P_{T_{E\setminus V}}(\cdot, x)$ (i.e., the measure $f\longmapsto
P_{T_{E\setminus V}}f(x)$) is concentrated on the boundary
$\partial V$ of $V$. Indeed, by (10.6) from \cite{Sh88} $
W_{T_{E\setminus V}}$ belongs to $E\setminus V$ a.s. on
$[T_{E\setminus V} <\infty]$. On the other hand we have
$T_{E\setminus V}>0$ $P^x$-a.s. and clearly $W_t(\omega)\in V$
provided that $t< T_{E\setminus V}(\omega)$. By the path
continuity of $\cw$ we conclude that $W_{T_{E\setminus V}}\in
\partial V$ $P^x$-a.s.
\end{rem}

Following \cite{God72}, a real-valued function $f$ defined on
an open set $V\subset E$ is called \emph{harmonic} on $V$, if it is
locally bounded, Borel measurable, finely continuous and there exists
$\rho>0$ such that
$$
f(x) = P_{T_{E\setminus B_r(x)}} f(x)
$$
for all $r<\rho$ whenever $\bar B_r(x) \subset V$; $\bar B_r(x)$
denotes the closed ball or radius $r$ centered at $x$.

We shall denote by $H^V:p\mathcal B(\partial V)\longrightarrow p \mathcal B(V)$
the kernel defined by
$$
H^V f := P_{T_{E\setminus V}}\bar f|_V\;,  \quad f\in p\mathcal
B(\partial V) ,
$$
where $\bar f$ is a Borel measurable extension of $f$ to $E$;
$H^Vf$ is well defined by Remark \ref{rem:hitting}.
 Hence
$$
H^Vf(x)  = E^x(f\circ W_{T_{E\setminus V}}\,;\, T_{E\setminus
V}<\infty)\;,  \quad x\in V.
$$
$H^V f$ is called the \emph{stochastic solution of the Dirichlet
problem} for $f$ (cf.\ \cite{God72}).

Recall that  (cf. \cite{Gr67}) an open set $V$ is called {\it
strongly regular} provided that for each $y\in\partial V$ there
exists a cone $K$ in $E$ with vertex $y$ such that $V\cap
K=\emptyset$; a cone in $E$ with vertex $y$ is the closed convex
hull of the set $\{y\}\cup\bar B_r(z)$ and $y\notin \bar
B_r(z)$.\\

By Corollary~1.2 and Remark~3.4 in
\cite{Gr67} it follows that:\\
\noindent $(5.1)  \quad$  if $V$ is strongly regular  and $f\in
\mathcal C(\partial V)$ is bounded, then $H^Vf$ is harmonic on $V$ and
$\lim\limits_{V\ni x\to y}H^Vf(x)=f(y)$ for all $y\in\partial
V$.\\

\noindent $(5.2)  \quad$ If $f\in \mathcal B(\partial V)$ is
bounded, then $H^Vf$ is harmonic on $V$ (see also Remark 3.4 in
\cite{Gr67} and page 453 in \cite{God72}).
Consequently, for every $f\in p\mathcal B(\partial V)$, $H^Vf$
is the sum of a series of positive harmonic functions on $V$.\\

\noindent {\bf Proof of $(5.2)$.} We may assume that $f\geq 0$. By
Theorem 3.6.4 in \cite{BeBo04} it follows that $H^Vf$ is an
excessive function with respect to the process on $V$ obtained  by
killing $\cw$ at the boundary of $V$. Therefore $H^Vf$ is finely
continuous on $V$ and $H^B H^Vf\leq H^Vf$ for all $B:= B_r(x)$,
$\overline{B}_r(x)\subset V$. Since $H^BH^V1(x)=H^V1(x)$ we
conclude that $H^BH^Vf(x)=H^Vf(x)$, hence $H^Vf$ is harmonic on
$V$. If $f\in p\cb(\partial V)$ then $H^Vf=\sum_n H^V f_n$, where
$(f_n)_n \subset bp\cb(\partial V)$ is such that $f=\sum_n f_n.$

\subsection*{Controlled convergence}
Let $f:\partial V\to\overline{\R}$, $V_0\subset V$, and
$h,k:V\to\overline{\R}$ be such that $k\geq 0$ and $h|_{V_0}, k|_{V_0}$ are real valued.
We say that $h$
\emph{converges to $f$ controlled by $k$ on $V_0$}, if the following
conditions hold: For every set $A\subset V_0$ and $y\in\partial
V\cap \bar A$ we have
\begin{enumerate}
\item[(c1)]
  If $\limsup\limits_{A\ni x\to y} k(x)<\infty$,
  then $f(y)\in\mathbb R$ and $f(y) = \lim\limits_{A\ni x\to y}h(x)$.
\item[(c2)]
  If $\lim\limits_{A\ni x\to y} k(x) =\infty$,
  then $\lim\limits_{A\ni x\to y}\frac{h(x)}{1+k(x)}=0$.
\end{enumerate}

\begin{rem} \label{rem.5.1}  
$(i)$ Following \cite{Co95} and \cite{Co98}, the controlled convergence
intends to offer a new method for setting and solving the
Dirichlet problem for general open sets and general boundary data.
In the above definition the function $f$ should be interpreted as
being the boundary data of the harmonic function $h$.
The function $k$ is called {\rm control function}, it is controlling
the convergence of the  solution $h$ to the given boundary data
$f$. If $\alpha> 0$ then $\alpha k$ and any majorant of $k$ are
also control functions.

$(ii)$ The case $k=0$, $V_0=V$,  corresponds to
the classical solution: $\lim\limits_{V\ni x\to y} h(x)=f(y)$ for
any boundary point $y$.

$(iii)$ In \cite{Co95} it was considered only the case $V_0=V$ for the controlled convergence.
It turns out that for the application we present here (see Theorem \ref{thm:t9}
below) we need to take into account an exceptional set $V\setminus V_0$.
\end{rem}

\noindent $(5.3)  \quad$  If $h_n$ converges to $f_n$ controlled
by $k$ on $V_0$ for each $n$ and $(\alpha_n)_n\subset \R$, $\alpha_n
\nearrow +\infty$, is such that $l:= \sum_n \alpha_n
|h_n|<\infty$, and $\sum_n h_n <\infty$ on $V_0$, then $\sum_n h_n$
converges to $\sum_n f_n$ controlled by $k+l$ on $V_0$  (cf. Proposition 1.7
in \cite{Co98}).

\begin{thm}\label{thm:t9} 
  Let $V\subset E$ be a strongly regular open set,
  $\lambda$ be a finite measure on $V$,  $\widehat \lambda$ be the
  measure on $\partial V$ defined by $\widehat \lambda := \lambda\circ H^V$,
and let  $f\in \mathcal L_+^1(\widehat \lambda)$.
Then there exist $g\in p\mathcal B(\partial V)$ and
a $\lambda$-zero set $M\subset V$  which is finely closed and $\lambda$-polar
with respect to the Brownian motion on $V$ (killed at the hitting time of $\partial V$),
such that
  $k:= H^Vg\in \mathcal L_+^1(\lambda)$ and
  $H^Vf$ converges to $f$ controlled by $k$ on  $V\setminus M$.
\end{thm}

\begin{pf}
Let $\mathcal{M} = \{ f \in \mathcal{L}_+^1 (\widehat \lambda): $
$\exists \; g \in p \mathcal{B} (\partial V)$ such that  $H^V f$
converges to $f$ controlled by $k = H^V g\in \cl^1(\lambda)$  on
$[k < \infty]\}$. Note that by $(5.1)$ the set of all positive
bounded continuous functions on $\partial V$ is a subset of
$\mathcal{M}$ (taking $k = 0$). Note also that the $\lambda$-zero
set $[k=\infty]$ is finely closed $\lambda$-polar because $k$ is a
$0$-excessive function with respect to  the Brownian motion on
$V$. The proof will be complete if we show that $\mathcal{M}$ is a
monotone class in $\mathcal{M}$.

Let $(f_n)_{n \geq 1} \subset \mathcal{M}$ be increasing to $f \in
\mathcal{L}_+^1 (\widehat \lambda)$. We show that $f \in
\mathcal{M}$. Let $h_n = H^V f_n$ and $h = H^V f$. Then $(h_n)_n$
increases to $h \in \mathcal{L}^1_+ (\lambda)$ and by hypothesis
$h_n$ converges to $f_n$ controlled by $k_n$  on $[k_n < \infty ]$ for all $n \geq 1$.
We may assume $\lambda (k_n) = 1$  for  all $n$. If
$$
k_0 : = \sum_n \frac{1}{2^n} k_n.
$$
then $h_n$ converges to $f_n$ controlled by $k_0$ on $[k_0<\infty]$ for all $n$. Let
$$
l: = \sum_{n \geq 1} n( h_{n + 1} - h_n) = \sum_{n \geq 1} (h - h_n).
$$
Since $\lambda (h_n) \nearrow \lambda (h) < \infty$, passing to a
subsequence, we may assume that $\sum_n (\lambda (h) - \lambda
(h_n)) < \infty$ and consequently $l = \mathcal{L}_+^1 (\lambda)$,
$l = H^V g$ with $g \in p \mathcal{B} (\partial V)$. By $(5.3)$ it
follows that $h$ converges to $f$ controlled by $k_0 + l$  on $[k_0
+ l < \infty]$, hence $f \in \mathcal{M}$.
\hfill $\square$
\end{pf}

\begin{rem} 
$(i)$ By $(5.2)$ the "solution" $H^Vf$ of the Dirichlet problem with boundary data
$f\in \cl^1_+(\widehat\lambda)$ from Theorem \ref{thm:t9}
is a sum of a series of positive harmonic functions on $V$.

$(ii)$ The result from Theorem \ref{thm:t9} holds in a more
general setting, e.g., for a path continuous Borel right process,
if $(5.1)$ holds.
\end{rem}

\section*{Appendix}
Let $(H, \langle , \rangle)$  be a separable real Hilbert space with norm
$|\cdot |$.  Let  $(E, \langle , \rangle_E)$ be another  Hilbert space
with norm $||\cdot ||$ such that $H\subset E$ continuously and densely  by a Hilbert-Schmidt map. Identifying $H$ with its  dual we have
$$
E' \subset H \subset E
$$
continuously and densely.
Let $\mu$ be a finitely additive measure on $H$ such that its Fourier transform
$\widehat\mu:H \longrightarrow \C$, defined by
$$
\widehat\mu(\xi):= \int_H e^{i\langle \xi , h\rangle} \mu (dh), \, \xi\in H,
$$
is continuous on $H$ and $\widehat\mu(0)=1$.
Then by the Bochner-Minlos Theorem (see, e.g., \cite{Ya89}) $\mu$ extends to
a probability measure on $(E, \cb(E))$ again denoted by $\mu$.\\

\noindent
{\bf Lemma A.1.} \label{lem:3.0}
{\it
Assume that apart from the Hilbert-Schmidt embedding $E' \subset H \subset E$
we have \mbox{another} such embedding
$$
E_1^{'} \subset H \subset E_1,
$$
i.e.,  $(E_1,\langle  \cdot , \cdot \rangle_{E_1})$ is a Hilbert
space with norm $||\cdot||_1:=\langle\cdot
,\cdot\rangle_{E_1}^\frac{1}{2}$ such that $H \subset E_1$
continuously and densely by a Hilbert-Schmidt embedding. Suppose
that there exists a linear subspace $K \subset E{'}\cap E'_1$ such
that $K$ separates the points both of $E_1$ and $E$ (i.e., for
each $x \in E\cup E_1$ such that $l(x)=0,$ for all $l\in K$, it
follows that $x=0$). Then there exists a Hilbert space
$(E_0,\langle\, , \, \rangle_{E_0})$ such that $H \subset E_0$
continuously and densely by a Hilbert-Schmidt map and both $E_0
\subset E$ and $E_0\subset E_1$ continuously. (Note that by
Kuratowski's theorem $E_0 \in \cb(E)\cap \cb(E_1)$.) }

\begin{pf}
Set $\pairing{}{h_1, h_2}{E_0}:=\pairing{}{h_1, h_2}{E} +
\pairing{}{h_1, h_2}{E_1}$, for all $h_1,h_2\in H$ with
corresponding norm $||\cdot||_{E_0}:={\pairing{}{\, ,
\, }{E_0}\!\!\!\!\!\!^\frac{1}{2}}$. Let $E_0$:=completion of $H$ with
respect to $||\cdot||_{E_0}$. Then clearly, $H \subset E_0$
continuously and densely by a Hilbert-Schmidt map.

\textit{Claim 1}. $E_0 \subset E$  continuously.

To prove the claim we have to show that
if $u_n\in H$, $n \in \N$, is an $||\cdot||_{E_0}$-Cauchy sequence
and at the same time an $||\cdot||_{E}$-zero sequence,
then it is also an $||\cdot||_{E_0}$-zero sequence.
But $u_n$, $n\in \N$, is also an $||\cdot||_{E_1}$-Cauchy sequence,
hence there exists $u \in E_1$ such that $\lim_{n\rightarrow \infty}||u_n-u||_{E_1}=0.$
It suffices to show that $u=0$.
To this end let $k \in K$.
Then $\pairing{E_1^{'}}{k , u}{E_1}=
\lim_{n\rightarrow \infty}\pairing{}{k , u_n}{H}=
\lim_{n\rightarrow \infty}\pairing{E'}{k , u_n}{E}=0$.

By assumption on $K$, it follows that $u=0$, and Claim 1 follows.

Likewise one proves:

\textit{Claim 2.} $E_0 \subset E_1$ continuously.

\hfill $\square$
\end{pf}

\noindent
{\bf Proposition A.2.}
{\it
Let $\{ e_n: n\in \N \}\subset E'$ be any orthonormal basis in $H$ separating
the points of $E$. For $n\in\N$ let $\widetilde P_n$ be defined by $(2.2)$
and $P_n:=\widetilde P_n \upharpoonright_H$.
Let  $\mu$ be a probability measure on $E$ coming from a cylinder measure on $H$,
i.e.,  $\mu$ is the image of a cylinder measure $\nu$ on $H$ under
the Hilbert-Schmidt embedding $H\subset E$,
and the Fourier transform $\widehat \nu$ of $\nu$ is continuous on $H$.
Then
\begin{center}
$\ds\lim_{n\rightarrow\infty}||z-P_nz||=0$ for $\mu$-a.e.
$z\in E$.
\end{center}
}

\begin{pf}
Let $\lambda_n\in(0,\infty)$ such that $\ds\sum_{n=1}^{\infty}\lambda_n<\infty$
and for $h_1,h_2\in H$ define
$$
\langle h_1,h_2\rangle_{E_1}:=\sum_{n=1}^{\infty}\lambda_n
\pairing{}{e_n, h_1}{H}\pairing{}{e_n, h_2}{H}
$$
with corresponding norm $||\cdot||_{E_1}:=\pairing{}{\cdot,
\cdot}{E_1}^{\frac{1}{2}}$. Let $E_1$ be the completion of $H$
with respect to $||\cdot||_{E_1}$. Then $H \subset E_1$
continuously and densely by a Hilbert-Schmidt map and hence we
have the Hilbert-Schmidt embeddings $E'_1 \subset H \subset E_1$.
 Furthermore
$\overline{e}_n:=\lambda_n^{-\frac{1}{2}}e_n$, $n \in \N$,
form an orthonormal basis of $E_1$ and for all $n \in \N$, $h \in
H$
$$
\leqno{(A.1)}\quad\quad \lambda_n \pairing{}{e_n, h}{H}=
\pairing{}{e_n, h}{E_1} \ ,
$$
hence
$$
\leqno{(A.2)}\quad\quad \pairing{E{'}}{e_n, h}{E}e_n
=\pairing{}{e_n, h}{H}e_n = \pairing{}{\overline{e}_n ,
h}{E_1}\overline{e}_n.
$$
Furthermore, for all $n \in \N$ by $(A.1)$
$$
h\longmapsto\pairing{}{e_n, h}{H}
$$
extends to a linear functional in $E'_{1}$ again denoted by
$e_n$. Hence $(A.1)$ implies by continuity that
$$
\leqno{(A.3)}\quad\quad \lambda_n \, {}_{E'_1}\!\!\pairing{}{e_n, z}{E_1}=
\pairing{}{e_n, z}{E_1}   \mbox{ for all } n\in \N, z\in E_1,
$$
in particular (since $\{\lambda_n^{-\frac{1}{2}}e_n\, : n\in
\N\}$ forms an ONB of $E_1$), $\{e_n :\,  n\in \N\}$ also
separates the points of $E_1$.
Hence   we can apply Lemma A.1
with  $K:=linspan\{e_n \; : \;n\in \N\}\subset E{'}$ (since $K$
also separates the points of $E$) to get the Hilbert space $E_0
\subset E \cap E_1$. Then the assertion of the proposition follows
from the following two claims.

\textit{Claim 1.} $\mu(E_0)=1.$

\textit{Claim 2.} $\ds\lim_{n\rightarrow\infty}{||P_nz-z||_{E}}=1,$
for all $z \in E_0.$

To prove Claim 1 we note that the cylinder measure on $H$ generating $\mu$,
mapped under the Hilbert-Schmidt embedding $H \subset E_0$ on
$E_0$, extends to a $\sigma$-additive probability measure on
$(E_0,\cb(E_0))$.
Clearly, because $H \subset E_0 \subset E$ continuously, we have
$\cb(E)\cap E_0=\cb(E_0)$, $E_0\in\cb(E)$, by Kuratowski's
theorem.
Hence it follows that this image measure coincides with $\mu$,
because the Fourier transforms coincide on $E{'}$ and $E' \subset E_0' \subset H
\subset E_0 \subset E$ continuously and densely. So, $\mu (E_0)=1$.

Now let us prove Claim 2.
By $(A.2)$ for all $h\in H$
$$
\leqno{(A.4)}\quad\quad  P_n h = \sum_{k=1}^{n}\pairing{}{\overline{e}_k, h}{E_1}\overline{e}_k.
$$

Let $z \in E_0$. Then there exists $h_l \in H$, $l \in \N$,
such that $\ds\lim_{l\rightarrow\infty}||z-h_l||_{E_0}=0$.
Hence, since both $E_0 \subset E$ and $E_0 \subset E_1$ continuously,
$$
\lim_{l\rightarrow\infty}||z-h_l||_{E}=0=
\lim_{l\rightarrow\infty}||z-h_l||_{E_1}.$$
Therefore, by $(A.4)$
$$
\leqno{(A.5)} \quad\quad\quad \widetilde P_n z= \ds\lim_{l\to\infty}P_n h_l=
\sum_{k=1}^{n}\pairing{}{\overline{e}_k,
z}{E_1}\overline{e}_k,\;\mbox{ for all } n\in \N.
$$
But the right hand side of $(A.5)$ converges to $z$,
since $\{\overline{e}_n : n\in \N\}$ is an orthonormal basis of
$(E_1,\pairing{}{\cdot, \cdot}{E_1})$, and Claim 2 is proved.
\hfill $\square$
\end{pf}

{\small \noindent {\bf Acknowledgement.}  Financial support by the
Deutsche Forschungsgemeinschaft, through Project GZ:436 RUM
113/23/0-1 and the CRC 701 is gratefully acknowledged. The first
named author acknowledges support from the Romanian Ministry of
Education, Research, Youth and Sport (PN II Program, CNCSIS code
ID 209/2007). The third named author would also like to thank the
I. Newton Institute in Cambridge for a very pleasant stay during
which part of this work was done. }












\end{document}